\documentclass[11pt]{article}

\usepackage{amsmath,amsfonts,amssymb,amsthm,mathtools,bbm, mathrsfs, geometry}

\DeclareMathAlphabet{\mathboondoxfrak}{U}{BOONDOX-frak}{m}{n}

\usepackage[all]{xy}

\usepackage{amssymb}

\newtheorem{thm}{Theorem}[section]

\newtheorem{dfn}[thm]{Definition}

\newtheorem{lem}[thm]{Lemma}

\newtheorem{prop}[thm]{Proposition}

\newtheorem{remark}[thm]{Remark}

\newtheorem{cor}[thm]{Corollary}

\newtheorem{ex}[thm]{Example}

\newtheorem{question}[thm]{Question}

\def\sq{{\scriptscriptstyle \square}}

\def\bq{\begin{question}}
\def\bt{\begin{thm}}
\def\bp{\begin{prop}}
\def\blem{\begin{lem}}
\def\bd{\begin{dfn}}
\def\br{\begin{remark}}
\def\bc{\begin{cor}}
\def\bex{\begin{ex}}
\def\beqs{\begin{eqnarray*}}
\def\beq{\begin{eqnarray}}
\def\bi{\begin{itemize}}

\def\eq{\end{question}}
\def\et{\end{thm}}
\def\ep{\end{prop}}
\def\elem{\end{lem}}
\def\ed{\end{dfn}}
\def\er{\end{remark}}
\def\ec{\end{cor}}
\def\eex{\end{ex}}
\def\eeqs{\end{eqnarray*}}
\def\eeq{\end{eqnarray}}
\def\ei{\end{itemize}}

\def\c{\cdot}

\def\ov{\overline}
\def\r{\rangle}
\def\l{\langle}

\def\H{{\cal H}}
\def\K{{\cal K}}
\def\B{{\cal{B}}} 
\def\fB{{\mathfrak B}}

\def\w*{w^*-w^*}

\def\ra{\rightarrow}

\def\vp{\varphi}
 
\def\bs{\backslash}

\def\E{{\cal E}}

\def\fB{{\mathboondoxfrak B}}
\def\lB{\lambda_{\mathcal B}}
\def\nuB{\nu_{\mathcal B}}
\def\nub{\nu_{\mathcal B}}
\def\lb{\lambda_{\mathcal B}}
\def\nuba{\nu_{{\mathcal B}^\alpha}}
\def\lba{\lambda_{{\mathcal B}^\alpha}}

\def\BS{{\bf S}} 
\def\cC{{\mathcal C}}
\def\lluc{\lambda_{LUC}}
\def\nuluc{\nu_{LUC}}

\def\nuap{\nu_{AP}}
\def\lwap{\lambda_{WAP}}
\def\nuwap{\nu_{WAP}}
\def\llmc{\lambda_{LMC}}
\def\nulmc{\nu_{LMC}}
\def\siglmc{\sigma_{LMC}}
\def\sigluc{\sigma_{LUC}}
\def\rholuc{\rho_{LUC}}
\def\sigwap{\sigma_{WAP}}

\def\nufb{\nu_{\fB}} 
\def\sigfb{\sigma_{\fB}}

\def\Bluc{\B^{LUC}}
\def\Bap{\B^{AP}}
\def\Bwap{\B^{WAP}}
\def\Blmc{\B^{LMC}}

\def\U{{\mathcal U}}
\def\V{{\mathcal V}}
\def\W{{\mathcal W}}
\def\K{{\mathcal K}}
\def\SK{{\cal SK}}
\def\BG{{_BG}}

\def\LO{{\cal LO}}
\def\RO{{\cal RO}}
\def\Nf{{{\cal N}_f}}
\def\vpnm{{\vp_{NM}}}

\def\SB{{\cal S}^\B}
\def\DB{{\Delta_\B}}
\def\db{{\delta_\B}}

\makeatletter
\def\moverlay{\mathpalette\mov@rlay}
\def\mov@rlay#1#2{\leavevmode\vtop{%
   \baselineskip\z@skip \lineskiplimit-\maxdimen
   \ialign{\hfil$\m@th#1##$\hfil\cr#2\crcr}}}
\newcommand{\charfusion}[3][\mathord]{
    #1{\ifx#1\mathop\vphantom{#2}\fi
        \mathpalette\mov@rlay{#2\cr#3}
      }
    \ifx#1\mathop\expandafter\displaylimits\fi}
\makeatother

\newcommand{\bigcupdot}{\charfusion[\mathop]{\bigcup}{\cdot}}

\def\S{{\cal S}}

\def\C{\mathbb{C}}
\def\R{\mathbb{R}}

\begin{document}

\title{Totally Disconnected Semigroup Compactifications of Topological Groups}
\author{Alexander Stephens  and  Ross Stokke    }
\date{}
\maketitle

\begin{abstract}{\small    We introduce the notion of an introverted Boolean algebra $\cal B$ of closed-and-open subsets of a topological group $G$,  show that the  associated Stone space $(\nu_\B G, \nu_\B)$ is a  totally disconnected semigroup compactification of $G$, and show  that every totally disconnected semigroup compactification of $G$ takes this form. We identify and study the universal totally disconnected semigroup compactification, the universal totally disconnected semitopological semigroup compactification and the universal totally disconnected group compactification of  $G$. Our main results are obtained independently of Gelfand theory and well-known properties of the (typically non-totally disconnected) universal compactifications $G^{LUC}$, $G^{WAP}$ and $G^{AP}$, though we do employ Gelfand theory  to clarify the relationship between these familiar universal compactifications and  their totally disconnected counterparts. 
\smallskip

}
\end{abstract} 

\section{Introduction } 

Totally disconnected locally compact topological groups provide an important class of examples  in abstract harmonic analysis, operator algebras, topological group theory,  and Lie group theory. A  sample of references in which the theory is developed and applied is \cite{Bur, Bur-Moz1, Bur-Moz2, Cap-Mon, Cap-Wil, Glo, Laz, Pal, Rau,  Rib-Zal, Sha-Wil, Willis,  Wil}.  Meanwhile, semigroup compactifications have played an essential role in the study of topological groups and semigroups, e.g., see \cite{Ber-Jun-Mil,  Rup}. When studying semigroup compactifications of totally disconnected groups, a natural aim is to remain within the class of totally disconnected spaces.  More generally, in light of the broad interest in totally disconnected structures, it is natural to investigate the totally disconnected semigroup compactifications of  any topological group. Two such totally disconnected semigroup compactifications were previously studied in \cite{Ili-Spr}: the universal totally disconnected group compactification and the ``idempotent compactification" associated with the Fourier--Stieltjes algebra of  a locally compact group.   Herein we will study the totally disconnected semigroup compactifications of a topological group $G$ by identifying them with  Stone compactifications $(\nub G, \nub)$ associated with ``introverted" Boolean algebras $\B$ of closed-and-open subsets of $G$.  The enveloping semigroup $E(X)$ of a topological dynamical system $(G, X)$  was introduced by Ellis in \cite{Ell}.  The totally disconnected semigroup compactifications with the joint continuity property, studied herein,   correspond exactly with the enveloping semigroups $E(X)$ of dynamical systems $(G, X)$ when $X$ is a totally disconnected compact Hausdorff space; see \cite[Section 10]{Gla}, noting that $E(X)$ is totally disconnected when $X$ is so.  For any group $G$, its profinite completion is a totally disconnected topological group compactification of $G$ \cite{Rib-Zal, Wil}, which we identify with a particular Stone compactification. Recent papers in which profinite completions of non-compact  groups are studied include \cite{BriEtAl1,BriEtAl2}. 
  
Let $X$ be a  \it (always Hausdorff) \rm topological space. Recall that $X$ is totally disconnected if all of its components are singleton sets, and $X$ is zero-dimensional if $\fB(X)$, the collection of closed-and-open --- hereafter called \it clopen \rm ---  subsets of $X$, is an open base for $X$. When $X$ is locally compact, it is totally disconnected if and only if it is zero-dimensional  \cite[Theorem 3.5]{HR}.  We will sometimes use the shorthand ``TD" in place of ``totally disconnected". By a \it compactification \rm  of $X$, we will mean a pair $(\alpha X, \alpha)$ where $\alpha X$ is a compact (Hausdorff) space and $\alpha: X \ra \alpha X$ is a dense-range continuous map; we will call $(\alpha X, \alpha)$ a \it topologist's compactification \rm when $\alpha$ is a homeomorphism of $X$ onto $\alpha (X)$.  When $(\alpha X, \alpha)$ and $(\gamma X, \gamma)$ are two compactifications of $X$, we write  $(\alpha X, \alpha) \preceq_X (\gamma X, \gamma)$ and say that $(\alpha X, \alpha)$ is a factor of $(\gamma X, \gamma)$ if  there is a continuous (necessarily surjective) map $h: \gamma X \ra \alpha X$ such that $h \circ \gamma = \alpha$; when $h$ is a homeomorphism, we call $h$ a compactification isomorphism, write  $(\alpha X, \alpha) \equiv_X (\gamma X, \gamma)$, and say that the two compactifications are equivalent.  The notation  $(\alpha X, \alpha) \precneqq_X (\gamma X, \gamma)$ is used  to indicate that $(\alpha X, \alpha)$ is a non-equivalent factor of  $(\gamma X, \gamma)$.  We let ${\cal K}(X)$ denote the set (up to equivalence) of all compactifications of $G$, and denote the subset of totally disconnected (equivalently, zero-dimensional)  compactifications by ${\cal K}_0(X)$.

We recall that  a semigroup and topological space $S$ is called a \it right (left) topological semigroup \rm if $ s \mapsto st$ ($s \mapsto  ts$) is continuous on $S$ for each $t \in S$; $S$ is a \it semitopological semigroup \rm if it is both a left and right topological semigroup; and $S$ is a \it topological semigroup \rm when $(s,t) \ra st$ is jointly continuous.

Unless stated otherwise, $G$ will always denote a (Hausdorff) topological  group.  A pair $(\alpha G, \alpha) \in {\cal K}(G)$ is called a \it right topological semigroup compactification of $G$,  \rm  or simply a \it semigroup compactification of $G$, \rm  if $\alpha G$ is a right topological semigroup  and $\alpha: G \ra \alpha G$ is a continuous  homomorphism such that $\alpha(G)$ is contained in the topological centre, 
$$Z_t (\alpha G) = \{ x \in \alpha G:  y \mapsto xy \text{ is continuous on } \alpha G\},$$
of $\alpha G$; note that $\alpha(e_G)$ is automatically the identity for $\alpha G$. One similarly defines \it left topological semigroup compactifications \rm  of $G$.   We let ${\cal SK}(G)$ and ${\cal SK}_0(G)$ respectively denote the sets (up to equivalence) of all semigroup compactifications of $G$ and all  totally disconnected semigroup compactifications of $G$; $(\SK(G), \preceq_G)$ and $(\SK_0(G), \preceq_G)$ are  complete lattices \cite[3.2.2, 3.3.3 and 3.3.4]{Ber-Jun-Mil}. When $(\alpha G, \alpha) \preceq_G (\gamma G, \gamma)$, the continuous intertwining map $h: \gamma G \ra \alpha G$ is automatically a surjective semigroup homomorphism  \cite[Proposition 3.1.6]{Ber-Jun-Mil}.

 For $f \in CB(G)$, the $C^*$-algebra of complex-valued functions on $G$ with the uniform norm $\| \c \|_\infty$,  and $s, t \in G$, let $f\c s(t) = f(st)$, $s \c f(t) = f(ts)$. A $C^*$-subalgebra $\S(G)$ of $CB(G)$ is \it left (right) $m$-introverted \rm if for $m$ in   $\Delta(\S(G))$ --- the Gelfand spectrum of $\S(G)$ --- $f \in \S(G)$ and $s, t \in G$, $f \c s, s \c f \in \S(G)$ and $m \c f \in \S(G)$ ($f \c m \in \S(G)$) where $m \c f(s) = m ( f\c s)$ ($f \c m(s) = m ( s \c f)$). When $\S(G)$ is left (right) introverted, $\Delta(\S(G))$ is a right (left) topological semigroup with respect to the product defined by $m\sq n(f) = m (n \c f)$ ($m \diamond n(f) = n (f \c m)$) for $m,n \in \Delta(\S(G))$ and $f \in \S(G)$; letting $\delta_\S(s) (f) = f(s)$ for $s \in G$, $f \in \S(G)$, $(\Delta(\S(G)), \delta_\S)$ is a right (left) topological ``Gelfand" semigroup compactification of $G$.

Every semigroup compactification of $G$ is, up to equivalence,  a Gelfand compactification, so every totally disconnected semigroup compactification is a Gelfand compactification. %We  study totally disconnected semigroup compactifications  as Gelfand compactifications in Section 5. 
Similarly, every topologist's compactification of a topological space $X$ is a Gelfand compactification. Nevertheless, there are advantages to having different constructions of  compactifications, with the suitability of each construction depending on the situation. 
For example, there are  many  constructions of the Stone-{\v C}ech compactification $\beta X$ of $X$ and, in the theory of semigroup compactifications, there are many constructions of the almost periodic/Bohr compactification $G^{AP}$ of $G$; advantages of each construction depends on the context.   In the case of zero-dimensional topological spaces $X$ (equivalently, totally disconnected spaces when $X$ is locally compact), the  Stone-compactification approach is especially apt \cite{Sto, Por-Woo}. Our primary goal  is to show how Stone's sublime construction provides a convenient approach to the study of totally disconnected semigroup compactifications of $G$.   Prior to Section 5, we do not employ function spaces, Gelfand theory or any known properties of  familiar semigroup compactifications such as $G^{LUC}$, $G^{WAP}$  and $G^{AP}$. %We hope the simplicity of our arguments demonstrate the advantages of this approach.   

In Section 2, to any Boolean subalgebra $\B$ of $\fB(X)$ we study the associated Stone compactification $(\nub X, \nub)$. Since we do not require that $\B$ be an open base for $X$, or even that  $X$ be zero-dimensional, the results in Section 2 seem, strictly speaking, to be new, though parallels are mostly known  when $\B$ is an open base for $X$.   
In subsequent sections, $G$ is assumed to be a topological group. In Section 3, we introduce (left/right) introverted Boolean subalgebras $\B$ of $\fB(G)$,  and show that 
the totally disconnected semigroup compactifications of $G$ are precisely, up to equivalence, of the form $(\nub G, \nub)$ for some left introverted Boolean subalgebra $\B$  of $\fB(G)$.  In Section 4 we identify Boolean subalgebras $\Blmc$, $\Bluc$, $\Bwap$ and $\Bap$ of $\fB(G)$ such that the corresponding Stone compactifications 
$(\nulmc G, \nulmc)$,  $(\nuluc G, \nuluc)$,  $(\nuwap G, \nuwap)$ and  $(\nuap G, \nuap)$ are the  universal totally disconnected semigroup compactification, the  universal totally disconnected semigroup compactification with the joint continuity property, the universal totally disconnected semitopological semigroup compactification and the universal totally disconnected group compactification of $G$. Beyond identifying and establishing the universal properties of these compactifications, we do not attempt to exhaustively detail their (expected) properties. As an application of our construction, we do however 
 identify $(\nuap G, \nuap)$ with the profinite completion of $G$.  In Section 5, we identify the totally disconnected semigroup compactifications $(\nub G,\nub)$ as Gelfand compactifications. In particular, we identify $(\nu_* G, \nu_*)$ with quotients in $\SK(G)$ of the familiar compactifications $(G^*, \delta_*)$, where $*$ is $LMC$, $LUC$, $WAP$ or $AP$.  

 % ; in this way $(\nuap G, \nuap)$ is identified in many instances, since profinite groups (and profinite completions of discrete groups)  are important in group theory, number theory and analysis. Perhaps more importantly, our $(\nuap G, \nuap)$ provides a new construction of the profinite completion of any topological group $G$. 

%Our approach in Section 4 is  independent of the familiar universal (not necessarily totally disconnected) semigroup compactifications $(G^{LMC}, \delta_{LMC})$, the universal semigroup compactification of $G$, $(G^{LUC}, \delta_{LUC})$, the ``left uniformly continuous"  universal semigroup compactification of $G$ with the joint continuity property, $(G^{WAP}, \delta_{WAP})$, the ``weakly almost periodic" universal semitopological semigroup compactification of $G$,  and $(G^{AP}, \delta_{AP})$, the  ``almost periodic/Bohr" universal group compactification of $G$. However, in Section 5  we identify $(\nu_* G, \nu_*)$ with quotients in $\SK(G)$ of the familiar compactifications $(G^*, \delta^*)$, where $*$ is $LMC$, $LUC$, $WAP$ or $AP$.  

\section{Totally disconnected compactifications as Stone spaces}

A readable account of the following ideas, based on the work of  M.H. Stone \cite{Sto}, is found in Sections 3.1, 3.2 and 4.7 of \cite{Por-Woo}.

Let $\B$ be a Boolean algebra, i.e., suppose that $(\B, \wedge, \vee, \leq, 0,1)$ is a complemented distributive lattice.  Canonical examples are the power set $(\mathcal{P}(X), \cap, \cup, \subseteq, \emptyset, X)$  of a set $X$, and the collection $(\fB(X), \cap, \cup,\subseteq, \emptyset, X)$ of all clopen subsets of a topological space $X$. 
 A  subset $\U$ of $\B$ is a \it filter \rm on $\B$  if it does not contain $0$, is closed upwards in $(\B, \leq)$, and is closed under the $\wedge$ operation;  $\U$ is called an \it ultrafilter  \rm if it is a maximal filter with respect to containment. Equivalently, a filter $\U$ is an ultrafilter when it is prime, meaning that $A \in \U$ or $B \in \U$ whenever $A \vee B \in \U$; equivalently, for every $B \in \B$,   $\U$ contains $B$ or its complement \cite[Proposition 3.1(r)]{Por-Woo}.

\bd \rm Letting $\BS(\B)$ denote the collection of all ultrafilters on $\B$, the sets
$$ \lB(B) = \{ \U \in \BS(B): B \in \U\} \quad (B \in \B)$$
form an open base for a topology on $\BS(\B)$; $\BS(\B)$ equipped with this topology is called the \it Stone space of $\B$ \rm  \cite[3.2 (b) and (c)]{Por-Woo}.  \ed 

\medskip 

\noindent {\bf The Stone Representation Theorem} \cite[Theorem 3.2(d)]{Por-Woo}
Let $\B$ be a Boolean algebra. Then: \bi 
\item[(i)] $\BS(\B)$ is a compact zero-dimensional Hausdorff space;
\item[(ii)] $\lb[\B] = \{\lb(B): B \in \B\} = \fB(\BS(\B))$, the collection of \it all \rm clopen subsets of $\BS(\B)$; and
\item[(iii)] $\lb: \B \ra \fB(\BS(\B))$ is a Boolean algebra isomorphism. 
\ei 
Thus, up to a Boolean algebra isomorphism, Boolean algebras are precisely the collections of all clopen subsets, $\fB(X)$, of (compact, zero-dimensional, Hausdorff) topological spaces $X$. 

\bd  \rm Let $X$ be a topological space and let $\B$ be a Boolean subalgebra of $\fB(X)$. Let $\nuB X = \BS(\B)$ be the Stone space of $\B$, and define
$$\nuB: X \ra \nuB X: x \ra \nuB(x) = \U_x \quad \text{where} \quad \U_x = \{B \in \B: x \in B\}.$$
\ed

 When $X$ is zero-dimensional and $\B \subseteq \fB(X)$ is an open base for $X$, it is well known that $(\nuB X, \nuB)$ is a zero-dimensional (equivalently, totally disconnected)  topologist's compactification of $X$ and  $(\beta_0 X, \beta_0):= (\nu_{\fB(X)} X, \nu_{\fB(X)})$ is the projective maximum in the set (up to equivalence) of all totally disconnected topologist's compactifications of $X$ \cite[Propositions 4.7 (b) and (c)]{Por-Woo}. Since  we do \it not \rm assume zero-dimensionality,   we now briefly develop some properties of $(\nuB X, \nuB)$  that we will need.  Parts (iv) and (v) of Proposition \ref{nuB X Prop} are almost certainly known in the special case that  $\cal A$ and $\B$ are  open bases for zero-dimensional $X$ and $(\alpha X, \alpha)$ is a topologist's compactification of $X$, but we were unable to find references. 
 
 \bp \label{nuB X Prop} Let $X$ be a topological space, $\B$ a Boolean subalgebra of $\fB(X)$. 
\bi \item[(i)] The pair $(\nuB X, \nuB)$ belongs to  $\K_0(X)$.
\item[(ii)] The compactification $(\nuB X, \nuB)$ is a topologist's compactification if and only if $\B$ is an open base for $X$ (and therefore $X$ is zero-dimensional).
\item[(iii)] For $B \in \B$, $\nuB(B) = \lB(B) \cap \nuB(X)$, $\nuB^{-1}(\lB(B)) = B$,  $cl_{\nuB X}(\nuB(B)) = \lB(B)$, and the map $\lB: \B \ra \fB(\nuB X): B \ra \lB(B)$ is a Boolean algebra isomorphism with inverse $D \mapsto \nuB^{-1}(D)$.
\item[(iv)] If $\cal A$ is another Boolean subalgebra of $\fB(X)$, then $(\nu_{\cal A} X, \nu_{\cal A}) \preceq_X (\nuB X, \nuB)$ if and only if ${\cal A} \subseteq \B$; in this case $h: \nub X \ra \nu_{\cal A} X: \U \mapsto \U \cap {\cal A}$ is a continuous surjection such that $h \circ \nub = \nu_{\cal A}$. Hence, $(\nu_{\cal A} X, \nu_{\cal A}) \equiv_X (\nuB X, \nuB)$ if and only if ${\cal A} = \B$. 
\item[(v)] If $(\alpha X, \alpha) \in \K_0(X)$, $\B^\alpha = \{ \alpha^{-1}(D): D \in \fB(\alpha X)\}$ is a Boolean subalgebra of $\fB(X)$ and $(\alpha X, \alpha) \equiv_X (\nu_{\B^\alpha} X, \nu_{\B^\alpha})$; $\B^\alpha$ is the unique Boolean subalgebra of $\fB(X)$ with this property.
\item[(vi)] Letting $\fB= \fB(X)$, $(\nu_\fB X, \nu_\fB)$ is the universal totally disconnected compactification of $X$, i.e, it is the projective maximum of  $\K_0(X)$. 
\ei 
\ep 

\begin{proof} (i) For $x \in X$, $x \in \nuB^{-1}(\lB(B))$ if and only if $B \in \nuB(x)$, equivalently $x \in B$; hence $\nuB^{-1}(\lB(B)) = B$. Since $\lB[\B]$ is an open base for $\nuB X$, $\nuB$ is continuous. Observing that $\nuB(B) = \lB(B) \cap \nuB(X)$, we have density of $\nuB(X)$ in $\nuB X$. 

\smallskip 

\noindent (ii) When $\B$ is an open base, it is easy to see that $\nuB$ is one-to-one and, since $\nuB(B) = \lB(B) \cap\nuB(X)$,  $\nuB$ is an open map of $X$ onto $\nuB(X)$. When  $\nuB: X \ra \nuB(X)$ is a homeomorphism,  the sets 
$ \nuB^{-1}(\lB(B) \cap \nuB(X)) = \nuB^{-1}(\lB(B)) = B$  $(B \in \B)$ compose an open base for $X$, i.e., $\B$ is an open base for $X$. 

\smallskip 

\noindent (iii) Let $B \in \B$. The first two identities were observed in (i), and the containment of $cl_{\nuB X}(\nuB(B))$ in the closed set $\lB(B)$ is obvious. If the open set $\lB(B) \bs cl_{\nuB X}(\nuB(B)) $ were nonempty, then it contains some $\nuB(x) = \U_x$ by density of $\nuB(X)$ in $\nuB X$. But then $B \in \U_x$, whence $x \in B$, a contradiction. The rest of (iii) follows immediately from the Stone representation theorem and the identity $\nuB^{-1}(\lB(B)) = B$ $(B \in \B)$. 

\smallskip 

\noindent (iv) When $\cal A$ is contained in $\B$, it is easy to check that $h: \nuB X \ra \nu_{\cal A} X: \U \mapsto \U \cap {\cal A}$ is well-defined, satisfies $h \circ \nuB = \nu_{\cal A}$, and for $A \in {\cal A}$, $h^{-1} (\lambda_{\cal A}(A)) = \lB(A)$, so $h$ is continuous. Hence, $(\nu_{\cal A} X, \nu_{\cal A}) \preceq_X (\nuB X, \nuB)$ in this case. Conversely, suppose that $h: \nuB X \ra \nu_{\cal A} X$ is a continuous map such that $h \circ \nuB = \nu_{\cal A}$. Let $A \in {\cal A}$. Then $h^{-1}(\lambda_{\cal A}(A))$ is a clopen subset of $\nuB X$, so part (ii) of the Stone representation theorem says that $h^{-1}(\lambda_{\cal A}(A)) = \lB(B)$ for some $B \in \B$. To see that $A = B$, suppose first that $x \in B$. Then $\nuB(x) \in \lB(B)$, so $\nu_{\cal A}(x) = h(\nuB(x))\in \lambda_{\cal A}(A)$; therefore, $x \in A$. On the other hand, if $x \in A$, then $h(\nuB(x))= \nu_{\cal A}(x) \in \lambda_{\cal A}(A)$, so $\nuB(x) \in h^{-1}(\lambda_{\cal A}(A)) = \lB(B)$; hence $x \in B$. Thus, $A = B \in \B$. 

\smallskip 

\noindent (v) Let $(\alpha X, \alpha)$ be a zero-dimensional compactification of $X$. Since $\alpha$ is continuous, $\B^\alpha = \{ \alpha^{-1}(D): D \in \fB(\alpha X)\}$ is a Boolean subalgebra of $\fB(X)$ and $\alpha^\sharp: \fB(\alpha X) \ra \B^\alpha: D \mapsto \alpha^{-1}(D)$ is a surjective Boolean algebra homomorphism. To see that $\alpha^\sharp$ is a Boolean algebra isomorphism, it now suffices to observe that by density of $\alpha(X)$ in $\alpha X$, $D \cap \alpha(X) \neq \emptyset$ for any nonempty set $D$ in $\fB(\alpha X)$, whence $\alpha^\sharp(D) = \alpha^{-1}(D) \neq \emptyset$.  Since $\alpha X$ is compact and filters satisfy the finite intersection property, one sees that $\BS(\fB(\alpha X)) = \{ \U_z: z \in \alpha X\}$, where $\U_z = \{ D \in \fB(\alpha X): z \in D\}$; because $\alpha X$ is zero-dimensional, $z \mapsto \U_z: \alpha X \mapsto \BS(\fB(\alpha X))$ is a bijection. Since $\alpha^\sharp$ is a Boolean algebra isomorphism, 
$\U_z \mapsto \alpha^\sharp(\U_z) = \{\alpha^{-1}(D): D \in \fB(\alpha X) \text{ and } z \in D\}$ 
is a bijection of $\BS(\fB(\alpha X))$ onto $\BS(\B^\alpha) = \nuba X$. Hence, 
$$h: \alpha X \ra \nuba X: z \mapsto \alpha^\sharp(\U_z) = \{\alpha^{-1}(D): D \in \fB(\alpha X) \text{ and } z \in D\}$$
is a bijection such that for $x \in X$, 
$$h(\alpha(x)) = \{\alpha^{-1}(D): D \in \fB(\alpha X) \text{ and } x \in \alpha^{-1}(D)\} = \{ B \in \B^\alpha: x \in B\} = \nuba(x).$$
Let $B \in \B^\alpha$, so $\lba (B)$ is a basic open subset of $\nuba X$; suppose that $B = \alpha^{-1}(D_0) = \alpha^\sharp(D_0)$ where $D_0 \in \fB(\alpha X)$. We will show that $h^{-1}(\lba(B)) = D_0$ (and conclude that $h$ is continuous). Suppose first that $z \in h^{-1}(\lba(B))$. Then $B \in h(z)$, so $B= \alpha^{-1}(D') = \alpha^\sharp(D')$ for some $D' \in \fB(\alpha X)$ such that $z \in D'$. Since $\alpha^\sharp$ is one-to-one, $D_0 = D'$, so $z \in D_0$. On the other hand, if $z \in D_0$, then $B = \alpha^{-1} (D_0) \in h(z)$; hence $h(z) \in \lba(B)$ and therefore $z\in h^{-1}(\lba(B))$. Thus, $h^{-1}(\lba(B)) = D_0$, as needed.  As a continuous map between compact Hausdorff spaces, $h$ is a homeomorphism. The uniqueness of $\B^\alpha$ is a consequence of part (iv). 

\smallskip 

\noindent (vi) This is an immediate consequence of parts (i), (iv) and (v). 
\end{proof} 

%Observe that parts (iv) and (v) of Proposition \ref{nuB X Prop}  combine to show that $\B \mapsto (\nub X, \nub)$ is an order isomorphism between the collection of all Boolean subalgebras of $\fB(X)$ and $\K_0(X)$. 

\section{Introverted Boolean algebras and semigroup compactifications} 

Throughout this section,   $G$ is a (Hausdorff) topological group. Define right and left group actions of $G$ on the collection of all subsets of $G$,  ${\cal P}(G)$, by putting 
$$B \c s = s^{-1}B \quad \text{ and } \quad s\c B = Bs^{-1} \qquad \text{ for } B \in {\cal P}(G), \ s \in G.$$
A Boolean subalgebra $\B$ of ${\cal P}(G)$ will be called \it $G$-invariant \rm if it is closed under these group actions; $\B$  will be called \it inversion invariant  \rm  if  $B^{-1} \in \B$ whenever $B \in \B$.  For example,  $\fB(G)$ is $G$-invariant and inversion invariant. For any collection ${\cal C}$ of subsets  of  $G$  and  $B \subseteq G$, let  
$${\cal C} \c B = \{ s \in G: B\c s \in \cC \} \quad \text{ and } \quad B \c \cC = \{ s \in G: s \c B \in \cC\}.$$

\bd \rm We will say that a $G$-invariant Boolean subalgebra $\B$ of $\fB(G)$ is \it left \rm  (respectively \it right\rm) \it introverted \rm  if 
$\V \c B \in \B$ (respectively $B \c \V \in \B$) for every $\V \in \nub G$ and $B \in \B$; we will say that $\B$ is \it introverted \rm if is both left and right introverted.
When $\B$ is left (respectively right) introverted, we define the \it left \rm  (respectively {\it right}) \it Arens product \rm on $\nub G$ by 
$$\U \sq \V = \{ B \in \B: \V \c B \in \U\} \quad \text{(respectively } \U \diamond \V = \{ B \in \B: B \c \U \in \V\}).$$ 
\ed  

We collect some elementary properties of these operations: 

\blem \label{Properties of operations Lemma} Let $\B$ be a $G$-invariant Boolean subalgebra of $\fB(G)$. \bi \item[(i)] For each $s \in G$, $B \mapsto B \c s$ and $B \mapsto s \c B$ are Boolean algebra isomorphisms of $\B$ onto itself.
\item[(ii)] For $\cC \subseteq \B$, $B \in \B$ and $s \in G$, $(\cC \c B) \c s = \cC\c (B \c s)$ and $s \c (B \c \cC) = (s \c B) \c \cC$. 
\item[(iii)] For $B \in \B$ and $s \in G$, $\nub(s) \c B = s \c B$ and $B \c  \nub(s) = B \c s$. 
\item[(iv)] For $\V \in \nub G$, $ B \mapsto \V \c B$  ($B \mapsto B \c \V$) defines a Boolean algebra homomorphism of $\B$ into ${\cal P}(G)$;  when $\B$ is left (right) introverted, we have a Boolean algebra homomorphism of $\B$ into itself. 
\item[(v)] When $\B$ is left (right) introverted, $(\U \sq \V) \c B = \U \c (\V \c B)$ ($B \c (\U \diamond \V) = (B \c \U) \c \V$) for $\U, \V \in \nub G$ and $B \in \B$. 
\ei 

\elem

\begin{proof}  Part (i) is clear and parts (ii) and (iii) are simple calculations. For example,  $t\in (\cC \c B) \c s = s^{-1} (\cC \c B)$ if and only if $st \in \cC \c B$ if and only if $(B\c s)\c t = B \c (st) \in \cC$  if and only if $t \in \cC \c ( B \c s)$.  Statement (iv) asserts that $\V \c \emptyset = \emptyset$, $\V \c G = G$, and for $B, C \in \B$, $\V \c (B \cup C)  = \V \c B \cup \V \c C$, $\V \c (B \cap C)  = \V \c B \cap \V \c C$ and $\V \c (G \bs B)  = G \bs (\V \c B)$, all of which follow from (i) and properties of ultrafilters. For example, ultrafilters are prime, so $s \in \V \c (B \cup C)$ if and only if $B\c s \cup C \c s = (B \cup C) \c s \in \V$ if and only if $B \c s \in \V$ or $C \c s \in \V$ if and only if $s \in \V \c  B \cup \V \c C$. Statement (v) follows from (ii): $s \in (\U \sq \V) \c B$ if and only if $B \c s \in \U \sq \V$ if and only if $(\V\c B) \c s = \V \c (B \c s) \in \U $ if and only if $s \in \U \c (\V \c B)$. 
\end{proof}  

\bt \label{First TD Semigroup Cpctn Thm} Let $G$ be a topological group. \bi 

\item[(i)] If $\B$ is a left introverted Boolean subalgebra of $\fB(G)$, then $(\nub G, \nub, \sq) \in \SK_0(G)$. 
\item[(ii)] If $(\alpha G, \alpha) \in \SK_0(G)$, then $\B^\alpha = \{ \alpha^{-1}(D): D \in \fB(\alpha G)\}$ is a left introverted  Boolean subalgebra of $\fB(G)$ and $(\alpha G, \alpha) \equiv_G (\nu_{\B^\alpha} G, \nu_{\B^\alpha})$;  $\B^\alpha$ is the unique Boolean subalgebra of $\fB(G)$ with this property.
\item[(iii)] If $\B$ is a  Boolean subalgebra of $\fB(G)$ and  $(\nub G, \nub)\in \SK_0(G)$ with respect to some product on $\nub G$, then $\B$ is left introverted and the product on $\nub G$ is the left Arens product.  
\ei
 \et 
 \begin{proof} (i) By Proposition \ref{nuB X Prop}, $(\nub G, \nub) \in \K_0(G)$. Let $\U,\V, \W\in \nub G$. That $\U \sq \V \in \nub G$ is a consequence of Lemma \ref{Properties of operations Lemma} (iv): $\V \c \emptyset = \emptyset \notin \U$, so $\emptyset \notin \U \sq \V$ and $\V \c G = G \in \U$, so $G \in \U \sq \V$. Suppose that $A, B \in \U \sq \V$. Then $\V \c A, \V \c B \in \U$, so $\V \c (A \cap B) = \V \c A \cap \V \c B \in \U$; therefore $A \cap B \in \U \sq \V$. If $C \in \B$ and $ C \supseteq B \in \U \sq \V$, then $\V \c C \supseteq \V \c B \in \U$, so $\V \c C \in \U$; therefore $C \in \U \sq \V$. For $B \in \B$, $\V \c B \in \U$ or $G\bs(\V \c B) = \V \c (G \bs B) \in \U$, so $B \in \U\sq \V$ or $G \bs B \in  \U \sq \V$. Hence, $\U \sq \V$ is a $\B$-ultrafilter. 
Lemma \ref{Properties of operations Lemma} (v) yields associativity of $\sq$: $B \in (\U \sq \V) \sq \W$ if and only if $\V \c (\W \c B) = (\V \sq \W) \c B \in \U$ if and only if $B \in \U \sq (\V \sq \W)$. 
To see that $\nub G$ is a right topological semigroup, suppose that $\U_i \ra \U$ in $\nub G$ and $\U \sq \V \in \lb(B)$ for some $B \in \B$. Then $\V \c B \in \U$, so $\U \in \lb(\V \c B)$. Therefore, $\U_i$ is eventually in $\lb(\V \c B)$ and $\U_i \sq \V$ is eventually in $\lb(B)$. Let $s, t \in G$. Then $\nub(s) \sq \nub(t) = \nub(st)$ by Lemma \ref{Properties of operations Lemma} (iii). Finally, to see that $\nub$ maps $G$ into $Z_t(\nub G)$, suppose that $\V_i \ra \V$ in $\nub G$ and $\nub (s) \sq \V \in \lb(B)$ for some $B \in \B$. Then $\V \c B \in \nub(s)$, so $B \c s \in \V$ meaning $\V \in \lb(B \c s)$. Hence, $\V_i$ is eventually in $\lb(B \c s)$ and therefore $\nub(s) \sq \V_i$ is eventually in $\lb(B)$, as required. 

\smallskip 

\noindent (ii)  By Proposition \ref{nuB X Prop} (v), we only need to show that $\B^\alpha$ is left introverted. Let $s \in G$, $B \in \B^\alpha$, say $B = \alpha^{-1}(D)$ where  $D \in \fB(\alpha G)$. Since $(\alpha G, \alpha)$ is a semigroup compactification of $G$, $r_{\alpha(s)}, \ell_{\alpha(s)}: z \mapsto \alpha(s) z, z \alpha(s)$ are homeomorphisms of $\alpha G$, with inverse maps $r_{\alpha(s^{-1})}, \ell_{\alpha(s^{-1})}$. Hence, $D \alpha(s^{-1}), \alpha(s^{-1}) D \in \fB(\alpha G)$ and therefore $s \c B = \alpha^{-1}(D) s^{-1} = \alpha^{-1}(D \alpha(s^{-1})) \in \B^\alpha$ and $B\c s = s^{-1} \alpha^{-1}(D) = \alpha^{-1}(\alpha(s^{-1})D ) \in \B^\alpha$. Thus, $\B^\alpha$ is $G$-invariant.   
Let $\V \in \nu_{B^\alpha} G$. From the proof of Proposition \ref{nuB X Prop} (v), $\V = h(z) = \{\alpha^{-1}(C): C \in \fB(\alpha G) \text{ and } z \in C\}$ for some $z \in \alpha G$. Since $r_z: \alpha G \ra \alpha G: x \mapsto xz$ is continuous, $r_z^{-1} (D) \in \fB(\alpha G)$ and therefore $\alpha^{-1}(r_z^{-1}(D)) \in \B^\alpha$. However, $s \in \V \c B$ if and only $B \c s =  \alpha^{-1}(\alpha(s^{-1})D ) \in \V = h(z)$. Recalling that $C \mapsto \alpha^{-1}(C)$ is one-to-one, this happens if and only if $z \in  \alpha(s^{-1})D $, which happens exactly when $\alpha (s) z = r_z(\alpha(s))\in D$,  equivalently, $s \in \alpha^{-1}(r_z^{-1}(D))$. Hence, $\V \c B \in \B^\alpha$. 

\smallskip 

\noindent (iii) By statement (ii) and Proposition \ref{nuB X Prop} (iv), this is obvious. 
 \end{proof} 

\br \rm \label{Right introverted version of Semigroup Cpctn Thm Remark}  When $\B$ is a right introverted Boolean subalgebra of $\fB(G)$, a symmetric argument shows that $(\nub G, \nub, \diamond)$ is a left topological semigroup compactification of $G$. Moreover, any left topological semigroup compactification  $(\alpha G, \alpha)$ of $G$  is equivalent to $(\nu_{\B^\alpha} G, \nu_{\B^\alpha})$, where $\B^\alpha = \{ \alpha^{-1}(D): D \in \fB(\alpha G)\}$ is a right introverted Boolean subalgebra of $\fB(G)$.  
\er 

We conclude this section with some  observations about the relationship between $\sq$ and $\diamond$. 

\bp \label{Top centre vs left and right Arens product Prop} Suppose that $\B$ is an introverted Boolean subalgebra of $\fB(G)$, $\U \in \nub G$. Then $\U \in Z_t(\nub G, \sq)$ if and only if $\U \sq \V = \U \diamond \V$ for every $\V \in \nub G$.  
\ep 

\begin{proof}  Suppose that $\U \in Z_t(\nub G, \sq)$ and let $\V \in \nub G$. Choose a net $(s_i)$ in $G$ such that $\nub (s_i) \ra \V$, so $\U \sq \nub(s_i)  \ra  \U \sq \V$ in $\nub G$. Let $B \in \B$ and suppose that $B \in \U \sq \V$. Then $\U \sq \V \in \lb(B)$, so $\U \sq \nub(s_i) \in \lb(B)$, eventually. Hence, $B \in \U \sq\nub(s_i)$, meaning $\nub(s_i) \c B = s_i \c B \in \U$ eventually; thus, $s_i \in B \c \U$, or $B\c \U \in \nub(s_i)$, and therefore $B \in \U \diamond \nub(s_i)$ eventually. Thus, $\U \diamond \nub(s_i)$ eventually belongs to the closed set $\lb(B)$, and therefore $\U \diamond \V \in \lb(B)$ since $\U \diamond \nub(s_i) \ra \U \diamond \V$. Hence $B \in \U \diamond \V$. If we begin by assuming that $B \in \U \diamond \V$, a symmetric argument shows that $B \in \U \sq \V$.  This proves the forward implication and the converse direction is trivial. 
\end{proof}

 We now introduce some topological notions to be used in Sections 4.1, 4.3, and the next proposition. 
Let $X$ be a topological space, $\B$ a Boolean subalgebra of $\fB(X)$. For $\V \in \nub X$, define $f_\V: \B \ra\{0,1\}$ by $f_\V(B) = 1$ if $B \in \V$, $f_\V(B) = 0$ if $B \notin \V$. For a subset  ${\bf A} $ of $\nub X$, let $\sigma(\B, {\bf A})$ be the weak topology on $\B$ induced by the functions $f_\V$ with $\V\in {\bf A}$. We are interested in the  \it weak \rm and \it point  topologies \rm on $\B$,  $\sigma_\B:= \sigma(\B, \nub X)$ and $\rho_\B := \sigma (\B, \nub (X))$.

\blem \label{Weak and point topologies lemma} Let $\B$ be a Boolean subalgebra of $\fB(X)$.  
\bi \item[(i)] The $\sigma_\B$- and $\rho_\B$-topologies on $\B$  are Hausdorff, zero-dimensional, and $\rho_\B$ is coarser than $\sigma_\B$.
\item[(ii)] Let $(B_i)_i$ be a net in $\B$, $B \in \B$. Then $B_i \ra B$ in $\sigma_\B$ if and only if for each $\V \in \nub X$, $B \in \V$ implies $B_i \in \V$ eventually, and $B \notin \V$ implies $B_i \notin \V$ eventually. 
\item[(iii)] Let $(B_i)_i$ be a net in $\B$, $B \in \B$. Then $B_i \ra B$ in $ \rho_\B$ if and only if for each $s \in X$, $s \in B $ implies $s \in B_i$ eventually, and $s \notin B$ implies $s \notin B_i$ eventually.
\item[(iv)] For any subset ${\bf A}$ of $\nu_\B X$,  the maps $\B \ra \B: B \mapsto X \bs B$, $\B \times \B \ra \B: (A, B) \mapsto A\cup B,\ A \cap B$  are $\sigma(\B, {\bf A})$-continuous. 
\item[(v)] If ${\cal C}$ is another Boolean subalgebra of $\fB(X)$ that contains $\B$, then $\sigma_\B = \sigma_{\cal C}{\large |}_\B$ and $\rho_\B = \rho_{\cal C}{\large |}_\B$.
\item[(vi)] If $G$ is a topological group and $\B$ is a $G$-invariant Boolean subalgebra of $\fB(G)$, then for each $t \in G$, $L_t, R_t: \B \ra \B: B \mapsto t \c B, B \c t$ are $\sigma_\B$-homeomorphic Boolean algebra isomorphisms.  
\item[(vii)] If $G$ is a topological group and $\B$ is an inversion-invariant Boolean subalgebra of $\fB(G)$, then $ \B \ra \B: B \mapsto B^{-1}$ is a  $\sigma_\B$-homeomorphic Boolean algebra isomorphism.  

\ei  
\elem

\begin{proof} (i) Since $\{0,1\}$ is zero-dimensional, $\sigma_\B$ and $\rho_\B$ are too, and $\rho_\B \subseteq \sigma_\B$ because $\nub(X) \subseteq \nub X$. If $B, C \in \B$ with $B \neq C$, say there is an element $s$ in $B \bs C$, then $f_{\nub(s)}^{-1}(\{1\})$ and $f_{\nub(s)}^{-1}(\{0\})$ are disjoint $\rho_\B$-open neighbourhoods of $B$ and $C$, respectively. Hence $\rho_\B$, and therefore $\sigma_\B$, is Hausdorff. 

\smallskip 

\noindent (ii) and (iii) For a subset ${\bf A}$ of $\nub X$, $B_i \ra B$ in $\sigma(\B, {\bf A})$ if and only if $f_\V(B_i) \ra f_\V(B)$ for each $\V \in {\bf A}$, so these statements (and the corresponding statement for $\sigma(\B, {\bf A})$) are clear. 

\smallskip 

\noindent   (iv) Suppose that $A_i \ra A$  and $B_i \ra  B$  in  $\sigma(\B, {\bf A})$. For $\V \in \nub X$ and $C \in \B$, exactly one of $C$ or $X \bs C$ belongs to $\V$, from which it follows that $X \bs B_i \ra X \bs B$ in $\sigma(\B, {\bf A})$. Let $\V \in {\bf A}$. If $A \cup B \in \V$, then $A \in \V$ or $\B \in \V$, so $A_i \in \V$ eventually or $B_i \in \V$ eventually; either way, $A_i \cup B_i \in \V$ eventually. Suppose that $A \cup B \notin \V$. Then $A \notin \V$ and $B \notin \V$, so there is an $i_0$ such that $A_i \notin \V$ and $B_i \notin \V$ for $i \succeq i_0$; since $\V$ is prime, $A_i \cup B_i \notin \V$ for $i \succeq i_0$.  Hence, $A_i \cup B_i \ra A \cup B$ in  $\sigma(\B, {\bf A})$. Similarly, one sees that $A_i \cap B_i \ra A \cap B$ in $\sigma(\B, {\bf A})$. 

\smallskip 

\noindent   (v) Suppose that $(B_i)$ is a net in $\B$ and $B \in \B$. By Proposition \ref{nuB X Prop} (iv), $h: \nu_{\cal C} X \ra \nub X: \U \mapsto \U \cap \B$ is a well-defined surjection. It  follows from (ii) that $B_i \ra B$ in   $\sigma_\B$ if and only if $B_i \ra B$ in $ \sigma_{\cal C}$. That $B_i \ra B$  in  $ \rho_\B$ if and only if $B_i \ra B$ in  $ \rho_{\cal C}$ is an obvious consequence of (iii). 

\smallskip 

\noindent (vi) Since $L_t, L_{t^{-1}}$ are Boolean algebra isomorphisms of $\B$, $L_{t^{-1}}(\V) \in \nub G$ exactly when $\V \in \nub G$. It follows that $t \c B_i \ra t \c B$ in  $\sigma_\B$ whenever $B_i \ra B$ in  $\sigma_\B$. 

\smallskip 

\noindent (vii) This is readily checked. 
\end{proof}

An introverted Boolean subalgebra $\B$ of $\fB(G)$ will be called \it Arens regular \rm if $\sq = \diamond$ on $\nub G$.

\bp \label{General Arens regular Prop}  Let $\B$ be a Boolean subalgebra of $\fB(G)$. The following statements are equivalent: \bi 
\item[(i)] $\B$ is left introverted and $(\nub G, \nub, \sq)$ is a semitopological semigroup; 
\item[(ii)] $\B$ is left introverted and for each $B \in \B$, $\nub G \ra (\B, \sigma_\B): \V \mapsto \V \c B$ is continuous; 
\item[(iii)]  $\B$ is introverted and Arens regular; 
\item[(iv)] $\B$ is right introverted and $(\nub G, \nub, \diamond)$ is a semitopological semigroup; 
\item[(v)] $\B$ is right introverted and for each $B \in \B$, $\nub G \ra (\B, \sigma_\B): \U \mapsto B \c \U$ is continuous. 
\ei 
\ep 

\begin{proof}  We will prove the equivalence of statements (i), (ii) and (iii); symmetric arguments establish the equivalence of statements (iii), (iv) and (v). Assume (i) holds. By Remark \ref{Right introverted version of Semigroup Cpctn Thm Remark}, there is a right introverted Boolean subalgebra $\cal C$ of $\fB(G)$ such that $(\nub G, \nub) \equiv_G (\nu_{\cal C} G, \nu_{\cal C})$. By Proposition \ref{nuB X Prop} (iv), $\B = {\cal C}$. Hence,  $\B$ is introverted and $\sq = \diamond$ by  Proposition \ref{Top centre vs left and right Arens product Prop}. Thus,  statement (i) implies statement (iii) and the converse is trivial. Suppose that $\V_i \ra \V$ in $\nub G$. To establish the equivalence of statements (i) and (ii),  we now show that 
 $\V_i \c B \ra \V \c B$ in  $\sigma_\B$   for each $B \in \B$
if and only if 
 $\U \sq \V_i \ra \U \sq \V$  in  $\nub G$  for each  $\U \in \nub G$:  
Assume the first of these conditions,  let $\U \in \nub G$, and suppose that $\U \sq \V \in \lb(B)$ for some $B \in \B$. Then $\V \c B  \in \U$, so $\V_i \c B \in \U$ eventually, whence $\U \sq \V_i \in \lb(B)$ eventually. Now assume that the second of these conditions holds, and let $B \in\B$. Let $\U \in \nub G$. If $\V \c B \in \U$, then $\U \sq \V \in \lb(B)$, whence $\U \sq \V_i \in \lb(B)$ eventually; therefore, $\V_i \c B \in \U$ eventually. If $\V \c B \notin \U$, then $\U \sq \V \in \nub G \bs \lb(B)$, so $\U \sq \V_i \in \nub G \bs \lb(B)$ eventually; therefore $\V_i \c B \notin \U$ eventually. 
\end{proof} 

The following is an immediate corollary of Theorem \ref{First TD Semigroup Cpctn Thm} and Proposition \ref{General Arens regular Prop}. 

\bc \label{Semitop semigp cpctn and Arens regularity Corollary}  Let $G$ be a topological group. Then $(\alpha G, \alpha)$ is a totally disconnected semitopological semigroup compactification of $G$ if and only if $(\alpha G, \alpha) \equiv_G (\nub G, \nub)$ for some (necessarily unique) introverted Arens regular Boolean subalgebra $\B$ of $\fB(G)$. 
 \ec

 Recall  that $(\alpha G, \alpha) \in \SK(G)$ is an \it involutive semigroup compactification of $G$ \rm if has a continuous involution $z\mapsto z^*$ --- meaning, $(z^*)^* = z$ and $(yz)^* = z^* y^*$ for all $y, z \in \alpha G$ --- such that $\alpha(s^{-1}) = \alpha(s)^*$ for all $s \in G$.

 \br \rm Suppose that $(\alpha G, \alpha)$, $(\gamma G, \gamma)$ are involutive semigroup compactifications of $G$ with $(\alpha G, \alpha) \preceq_G (\gamma G, \gamma)$; say $h: \gamma G \ra \alpha G$ is a continuous (homomorphic surjective) map such that $h \circ \gamma = \alpha$. One can readily check that $h$ automatically satisfies $h(x^*) = h(x)^*$ for $x \in \gamma G$.  
 \er

\bp \label{Involutive semigroup compactifications Prop}  Let $(\alpha  G, \alpha)$ be a totally disconnected  involutive semigroup compactification of $G$. Then $(\alpha G, \alpha)$ is a semitopological semigroup compactification of $G$ and $(\alpha G, \alpha)\equiv_G (\nub G, \nub)$ for a unique introverted,  Arens regular, inversion-invariant Boolean subalgebra $\B$ of $\fB(G)$. 
 \ep 
 
 \begin{proof}  Take $x \in \alpha G$. Since $\alpha G$ is right topological, $\alpha G \ra \alpha G: y\mapsto y^* \mapsto y^* x^* \mapsto (y^* x^*)^* = xy$ is continuous. Hence, $(\alpha G, \alpha)$ is semitopological, so $(\alpha G, \alpha) \equiv_G (\nu_{\B^\alpha} G, \nu_{\B^\alpha})$ where $\B^\alpha = \{ \alpha^{-1}(D): D \in \fB(\alpha G)\}$  is an introverted Arens regular Boolean subalgebra of $\fB(G)$ by Theorem \ref{First TD Semigroup Cpctn Thm} and Corollary \ref{Semitop semigp cpctn and Arens regularity Corollary}. Let $D \in \fB(\alpha G)$. Since $x \mapsto x^*$ is a homeomorphism on $\alpha G$, $D^* = \{x^*: x \in D\} \in \fB(\alpha G)$ and $\alpha^{-1}(D^*) = \alpha^{-1}(D)^{-1}$, so $\B^\alpha$ is inversion invariant.  
 \end{proof} 

When $(\nub G, \nub)$  is an involutive semigroup, we now describe the involution. 

\bp \label{Involution Prop} Let $\B$ be an inversion-invariant left introverted Boolean subalgebra of $\fB(G)$.  \bi 
\item[(i)] For $\U \in \nub G$, $\U^*: = \{B^{-1} : B \in \U\} \in \nub G$; $(\U^*)^* = \U$; $\U \mapsto \U^*$ is a homeomorphism of $\nub G$; and for each $s \in G$, $\nub(s^{-1}) = \nub(s)^*$. Moreover,  $\B$ is introverted and $(\U \sq \V)^* = \V^* \diamond \U^*$ for $\U, \V \in \nub G$.   
\item[(ii)]  Hence, $\B$ is Arens regular if and only if $(\nub G, \nub, \sq)$ --- equivalently $(\nub G, \nub, \diamond)$ --- is an involutive semigroup compactification of $G$. 
\ei 
\ep 

\begin{proof}  (i) Since $B \mapsto B^{-1}$ is a Boolean algebra isomorphism on $\B$,  $\U^* \in \nub G$ for $\U \in \nub G$; clearly  $(\U^*)^*=\U$ and $\nub(s^{-1}) = \nub (s)^*$ for $s \in G$. To see that $\U \mapsto \U^*$ is continuous (and therefore a homeomorphism), observe that the  pre-image of $\lb(B)$ under this map is $\lb(B^{-1})$.  Let $\U, \V \in \nub G$. For $B \in \B$,  observe that  $\V \c B^{-1} = (B \c \V^*)^{-1}$. Since $\B$ is inversion invariant and left introverted, we deduce that $\B$ is also right introverted. Moreover, for $B \in \B$,   $B \in (\U \sq \V)^*$ if and only if $(B \c \V^*)^{-1} = \V \c B^{-1} \in \U$; equivalently, $B\in \V^* \diamond \U^*$.  

\smallskip 

\noindent (ii)  If $\B$ is Arens regular, we obtain $(\U \sq \V)^* = \V^* \diamond \U^* = \V^* \sq \, \U^*$. If  $(\nub G, \sq)$ is involutive, $\U \sq \V = (\U \sq \V)^{**}  = (\V^* \sq \, \U^*)^* = \U^{**} \diamond \V^{**} = \U \diamond \V$. %; if $(\nub G, \diamond)$ is involutive, $\U \diamond \V = (\U \diamond \V)^{**}  = (\V^* \diamond \U^*)^* = (\U \sq \V)^{**} = \U \sq \V$.
 \end{proof} 

 In Section 4 we will  observe that there exist non-Arens regular  introverted  inversion-invariant  Boolean subalgebras of $\fB(G)$, and  we will identify the largest Arens regular (inversion-invariant) Boolean subalgebra of $\fB(G)$. 
 
\section{Universal compactifications} 

For a topological group $G$, we will identify several totally disconnected semigroup compactifications as 
% universal totally disconnected semigroup compactification $(\nulmc G, \nulmc)$ of $G$, the universal totally disconnected semigroup compactification with the joint continuity property $(\nuluc G, \nuluc)$ of $G$, the universal totally disconnected semitopological semigroup compactification $(\nuwap G, \nuwap)$  of $G$, and the universal totally disconnected semigroup compactification $(\nuap G, \nuap)$ of $G$; each is identified as
  Stone spaces of  associated left introverted Boolean subalgebras of $\fB(G)$. All proofs are independent of function spaces and Gelfand theory, and of previously known universal properties of other (not necessarily totally disconnected) semigroup compactifications. Broad ideas from the non-totally disconnected theory, e.g. see \cite{Ber-Jun-Mil}, are sometimes employed.

\subsection{The universal TD semigroup compactification}

Let $G$ be a topological group. Writing $\nufb$, $\sigfb$  in place of $\nu_{\fB(G)}$, $\sigma_{\fB(G)}$ (and so on), let 
$$\Blmc := \{ B \in \fB(G):  s \mapsto B \c s: G \ra (\fB(G), \sigfb) \text{ is continuous}\}.$$
Using Lemma \ref{Weak and point topologies lemma} (iv) and (vi), one sees that $\Blmc$ is a $G$-invariant Boolean subalgebra of $\fB(G)$.

\bigskip 

\noindent {\bf Notation}: We will write $(\nulmc G, \nulmc)$, $\llmc(B)$, $\siglmc$, etc., in place of $(\nu_{\Blmc} G, \nu_{\Blmc})$, $\lambda_{\Blmc}(B)$, $\sigma_{\Blmc}$, etc.  Without further comment, we will  employ similar notation after later introducing the Boolean subalgebras $\Bluc$, $\Bwap$ and $\Bap$ of $\fB(G)$.

\bt \label{Universal semigp cpctn Thm}  Let $G$ be a topological group. Then $\Blmc$ is the largest left introverted Boolean subalgebra of $\fB(G)$ and $(\nulmc G, \nulmc)$  is  the universal totally disconnected semigroup compactification of $G$; i.e., $(\nulmc G, \nulmc)$ is the projective maximum of $\SK_0(G)$. 
 \et 
 
 \begin{proof}   Let $\V \in \nulmc G$, $B \in \Blmc$. Let $(t_i)$ be a net in $\V \c B$ such that $t_i \ra t$ in $G$. Then each $ B \c t_i \in \V$ and, by Lemma \ref{Weak and point topologies lemma} (v),  $B \c t_i \ra B \c t$ in $\siglmc$; hence $B \c t \in \V$, meaning that $t \in \V \c B$. Thus, $\V \c B$ is closed. Since $G \bs (\V \c B) = \V \c (G \bs B)$ is also closed, $\V \c B \in \fB(G)$.  To see that $\V \c B\in \Blmc$, suppose that $s_i \ra s$ in $G$ and let $\U \in \nufb G$. Define $\U \tilde{\sq} \V = \{ A \in \Blmc: \V \c A \in \U\}$. By Lemma \ref{Properties of operations Lemma} (iv), $A \mapsto \V \c A: \Blmc \ra \fB(G)$ is a Boolean algebra homomorphism, so the argument found in the first paragraph of the proof of Theorem \ref{First TD Semigroup Cpctn Thm} (i) showing there that $\U \sq \V$ is a $\B$-ultrafilter, shows here that $\U \tilde{\sq} \V$ is a $\Blmc$-ultrafilter, i.e., $\U \tilde{\sq} \V \in \nulmc G$. If $(\V \c B) \c s = \V \c (B \c s) \in \U$, then $B \c s \in \U \tilde{\sq} \V$ and therefore, since $B \c s_i \ra B \c s$ in $\siglmc$ (again by Lemma \ref{Weak and point topologies lemma} (v)),  $B \c s_i \in \U \tilde{\sq} \V$, eventually; hence $(\V \c B)\c s_i = \V\c (B \c s_i) \in \U$, eventually.   If $(\V \c B) \c s = \V \c (B \c s) \notin \U$, then $B \c s \notin \U \tilde{\sq} \V$, so  $B \c s_i \notin \U \tilde{\sq} \V$, eventually; hence $(\V \c B)\c s_i = \V\c (B \c s_i) \notin \U$, eventually. Thus, $(\V \c B) \c s_i \ra (\V \c B) \c s$ in  $\sigfb$, and we conclude that $\V \c B \in \Blmc$. Thus, $\Blmc$ is a left introverted Boolean subalgebra of $\fB(G)$ and, therefore, $(\nulmc G, \nulmc) \in \SK_0(G)$. 

Let $(\alpha G, \alpha) \in \SK_0(G)$. By Theorem \ref{First TD Semigroup Cpctn Thm}, $(\alpha G, \alpha) \equiv_G (\nub G, \nub)$ for some left introverted Boolean subalgebra $\B$ of $\fB(G)$. Let $B \in \B$, suppose that $s_i \ra s $ in $G$, and take $\V \in \nub G$. If $B \c s \in \V$, then $s \in \V \c B$, meaning $\V \c B \in \nub(s)$, and therefore $B \in \nub(s) \sq \V$; hence $\nub(s) \sq \V \in \lb(B)$. Since $\nub(s_i) \sq \V \ra \nub(s) \sq \V$ in $\nub G$, $\nub(s_i) \sq \V \in \lb(B)$ eventually, and therefore $B \c s_i \in \V$ eventually. If $B \c s \notin \V$, then $\nub(s) \sq \V$ belongs to the clopen set $\nub G \bs \lb(B)$, whence $\nub(s_i) \sq \V \in \nub G \bs \lb(B)$ eventually; therefore, $B \c s_i \notin \V$ eventually. Thus, $B \c s_i \ra B \c s $ in  $\sigma_\B$, so $B \c s_i \ra B \c s$ in $\sigfb$. Hence, $\B$ is contained in $\Blmc$ and therefore $(\nub G, \nub) \preceq_G (\nulmc G, \nulmc)$.% by Proposition \ref{nuB X Prop} (iv).
 \end{proof}

\subsection{Universality with respect to the joint continuity property} 

For a subset $E$ of $G$, let $$_E G= \{ s \in G: E \c s = E \} \quad \text{ and } \quad G_E= \{ s \in G: s \c E = E \}$$
be the right and left stabilizer subgroups of $G$.  Let 
$$\Bluc = \{ B \in \fB(G): \BG \text{ is open} \}.$$
%Though the results in this section are independent of function spaces, to provide context we observe below that $B \in \Bluc$ exactly when the characteristic function $1_B$ belongs to $LUC(G)$, the $C^*$-algebra of uniformly continuous functions on $G$. (Recall that $f \in LUC(G)$ if $s \mapsto f \c s : G \ra (CB(G), \|\c \|_\infty)$ is continuous, where $f \c s(t) = f(st)$ for $s,t \in G$.) 

\blem \label{BLUC Lemma} (a) Let $B \in \fB(G)$. The following statements are equivalent: 
\bi 
\item[(i)] $B \in \Bluc$;
%\item[(ii)] $_BG \in \Bluc$;
\item[(ii)] $G \ra \fB(G) : s \mapsto B\c  s$ is continuous (at $e_G$) with respect to the discrete topology on $\fB(G)$;
\item[(iii)] for each ${\cal A} \subseteq \fB(G)$, ${\cal A} \c B = \{ s\in G: B \c s \in {\cal A}\}$ is open (equivalently closed, equivalently clopen) in $G$; 
\item[(iv)]  for each $ C \in  \fB(G)$, $\{ C \} \c B = \{ s\in G: B \c s = C\}$ is open in $G$.
\ei
(b) The collection of sets $\Bluc$ is a left introverted Boolean subalgebra of $\fB(G)$. 
\elem
\begin{proof} (a) %Since $H = \BG$ is a subgroup of $G$, $_HG= H $ is open if and only if $H =\BG$ is open, i.e., statements (i) and (ii) are equivalent. 
Suppose that $B \in \Bluc$ and $s_i \ra e_G$ in $G$. Since $e_G$ belongs to the open subgroup $_BG$, $s_i \in \BG$, and therefore $B \c s_i =B\c e_G$, eventually. If $s_i \ra s$, then $s_i s^{-1} \ra e_G$ and therefore $B \c (s_i s^{-1}) =B$, whence $B \c  s_i = B\c s$, eventually. Hence,  (i) implies (ii). Since ${\cal A} \c B$ is the pre-image of  $\cal A$ under the map $s \mapsto B \c s$, statements (ii) and (iii) are equivalent. It is clear that (iii) implies (iv) and, since $\{ B \} \c B = \BG$, (iv) implies (i). 
\smallskip 

\noindent (b) Observe that ${_\emptyset G}=G$ and ${_GG}= G$, so $\emptyset, G \in \Bluc$. Suppose that $A, B \in \Bluc$. Then $_A G\cap \BG$ is an open subgroup of both $_{A\cap B}G$ and $_{A\cup B}G$, so $A \cap B, A\cup B \in \Bluc$; also $_{G \bs B}G=\BG$, so $G\bs B \in \Bluc$. Hence, $\Bluc$ is a Boolean subalgebra of $\fB(G)$.  For $t \in G$, observe that $_{t\c B} G= \BG$ and $_{B \c t}G = t^{-1}\BG t$; therefore $\Bluc$ is $G$-invariant. For $\V \in \nuluc G$, $B \in \Bluc$ and $t \in G$, $\V \c B \in \fB(G)$ by part (a), and $(\V \c B) \c t = \V \c (B \c t)$, so $\BG \subseteq {_{\V \c B}G}$; hence $\V \c B \in \Bluc$. 
\end{proof} 

%We will write $(\nuluc G, \nuluc)$ in place of $(\nu_{\Bluc} G, \nu_{\Bluc})$ --- a totally disconnected semigroup compactification of $G$ by Lemma \ref{BLUC Lemma} --- and $\lluc$ in place of $\lambda_{\Bluc}$. 
 A semigroup compactification $(\alpha G, \alpha)$ has the \it joint continuity property \rm if the map $G \times \alpha G \ra \alpha G: (s, z) \mapsto \alpha(s) z$ is jointly continuous \cite[Section 4.4]{Ber-Jun-Mil}. 
 
 \br \label{Cpctn Ppts Prsvd By Qtnts Remark}  \rm Suppose that $(\gamma G, \gamma) \in \SK(G)$ has the  joint continuity property/is a semitopological semigroup/is a topological group. Then any quotient of $(\gamma G, \gamma)$ in $\SK(G)$ has the same property: (Though this is likely known, we do not have a reference.) Suppose that $(\alpha G, \alpha)\preceq_G (\gamma G, \gamma)$ in $\SK(G)$, $h: \gamma G \ra \alpha G$ the associated compactification homomorphism. Assume first that $(\gamma G, \gamma)$ has the joint continuity property.  To see that  $(\alpha G, \alpha)$ has the joint continuity property, assuming that $(s_i , z_i) \ra (s, z)$ in $G \times \alpha G$, it suffices to show that $\alpha(s_{i_j})z_{i_j} \ra \alpha(s) z$ for some subnet $((s_{i_j}, z_{i_j}))_j$ of $((s_i, z_i))_i$. For each $i$, take $x_i$ in $\gamma G$ such that $h(x_i) = z_i$; take  $(x_{i_j})_j$ to be a subnet of $(x_i)_i$ and $x \in \gamma G$ such that $\lim x_{i_j} = x$. Then $\gamma(s_{i_j}) x_{i_j} \ra \gamma(s) x$ in $\gamma G$, so  $\alpha(s_{i_j}) z_{i_j} = h(\gamma(s_{i_j}) x_{i_j}) \ra h(\gamma(s) x) =  \alpha(s) z$, as needed.   A similar argument shows that $\alpha G$ is semitopological when $\gamma G$ is so. If $\gamma G$ is a topological group, then $\alpha G$ is a semitopological compact group, therefore a topological group by \cite[Corollary 1.4.5]{Ber-Jun-Mil}. 
 \er 

\bt \label{Universality of nuLUC Thm} Let  $\B$ be a Boolean subalgebra of $\fB(G)$. The following statements are equivalent: 
\bi \item[(i)]  $(\nub G, \nub)$ is a  semigroup compactification of $G$ with the joint continuity property;
\item[(ii)]  $\B$ is  left introverted and contained in $\Bluc$.
\ei  In particular,  $(\nuluc G, \nuluc)$ is universal among all totally disconnected semigroup compactifications of $G$ with the joint continuity property. 
\et 

\begin{proof}   Suppose that $(s_i, \V_i) \ra (s, \V)$ in $G \times \nuluc G$ and let $B \in \Bluc$ be such that $\nuluc(s) \sq \V \in \lluc(B)$; equivalently, $\V \in \lb(B \c s)$. Take $i_0$ such that $B \c s_i = B\c s$ for $i \succeq i_0$, and take $i_1 \succeq i_0$ such that $\V_i \in \lluc(B \c s)$ for $i \succeq i_1$. Then for  $i\succeq i_1$,   $B \c s_i = B \c s \in \V_i$, and therefore %whence $s_i \in \V_i \c B$, or $\V_i \c B \in \nuluc(s_i)$. Thus, $B \in \nuluc(s_i) \sq \V_i$, and therefore
 $\nuluc(s_i) \sq \V_i \in \lluc(B)$. Hence, $(\nuluc G, \nuluc)$ has the joint continuity property. 

Assume that    $(\nub G, \nub)$ is a  semigroup compactification of $G$ with the joint continuity property. By Theorem \ref{First TD Semigroup Cpctn Thm}, $\B$ is left introverted.  Suppose towards a contradiction that $B \notin \Bluc$ for some $B \in \B$. Then by Lemma \ref{BLUC Lemma}  there is a net $(s_i)$ converging to $e_G$ in $G$ such that $B \c s_i $ does not equal $B \c e_G = B$ eventually; passing to a subnet if needed, we can assume that $B \c s_i \neq B$ for all $i$. For each $i$, either $B \c s_i \bs B \neq \emptyset$ or $B \bs B \c s_i \neq \emptyset$, so there is some $t_i \in B$ such that either 
\beq  \label{Universality of nuLUC Thm Eqn} s_i^{-1} t_i \notin B \qquad \text{or} \qquad s_it_i \notin B.
\eeq
For each $i$, $\nub(t_i) \in \lb(B)$, a closed subset of the compact space $\nub G$. Passing to a subnet again, if necessary,  $\nub(t_i) \ra \V$ in $\nub G$ for some $\V \in \lb(B)$. Using (\ref{Universality of nuLUC Thm Eqn}) and passing yet again to a subnet if required, we can assume that either $s_i^{-1}t_i \notin B$ for each $i$, or $s_it_i \notin B$ for each $i$. Suppose first that  $s_i^{-1}t_i \notin B$ for each $i$. Since $(s_i^{-1}, \nub(t_i)) \ra (e_G, \V)$ in $G \times \nub G$ and $(\nub G, \nub)$ has the joint continuity property, $\nub(s_i^{-1} t_i) = \nub(s_i^{-1}) \sq \nub(t_i) \ra \nub(e_G) \sq  \V = \V$ in $\nub G$.  But $\V \in \lb(B)$, so for large $i$ we must have $\nub(s_i^{-1} t_i) \in \lb(B)$, meaning $s_i^{-1}t_i \in B$, a contradiction. One similarly obtains a contradiction if $s_it_i \notin B$ for each $i$. We conclude that $\B$ is contained in $\Bluc$, and therefore $(\nub G, \nub) \preceq_G (\nuluc G, \nuluc)$.   
By Theorem \ref{First TD Semigroup Cpctn Thm} we conclude that  $(\alpha G, \alpha)\preceq_G (\nuluc G, \nuluc)$ for any $(\alpha G, \alpha) \in \SK_0(G)$ with the joint continuity property. We have established the last statement of the theorem and that (i) implies (ii);  (ii) implies (i) is a consequence of Proposition \ref{nuB X Prop} (iv) and Remark \ref{Cpctn Ppts Prsvd By Qtnts Remark}.
\end{proof}  

When $G$ is a locally compact group, every semigroup compactification of $G$ has the joint continuity property by Ellis's joint continuity theorem, (which is contained in \cite[Theorem 1.4.2]{Ber-Jun-Mil}). We therefore have the following immediate consequence of Theorems \ref{Universal semigp cpctn Thm}, \ref{Universality of nuLUC Thm} and Proposition \ref{nuB X Prop} (iv). 

\bc \label{nuLUC universal for lcgs Corollary}  Let $G$ be a locally compact group. Then $(\nuluc G, \nuluc) \equiv_G (\nulmc G, \nulmc)$, the universal totally disconnected semigroup compactification of $G$, and $\Bluc=\Blmc$,  the largest left introverted Boolean subalgebra of $\fB(G)$. 
\ec

\br \rm   \label{Non-Arens regular example remark}   Symmetric arguments show  that $\B^{RUC} = \{ B \in \fB(G): G_B \text{ is open} \}$ is a right introverted Boolean subalgebra of $\fB(G)$ such that  the obvious  versions of  \ref{Universality of nuLUC Thm} and  \ref{nuLUC universal for lcgs Corollary}  hold for the left topological semigroup compactification $(\nu_{RUC} G, \nu_{RUC})$.

   %  \smallskip
 
 % \noindent (b) The question of whether there is an  example of a topological group for which $\Bluc \subsetneq \Blmc$ is left as an open problem. 

 %  \smallskip

 % \noindent  (c) Being compact and totally disconnected is productive and closed hereditary, and is therefore preserved by the formation of subdirect products. By \cite[Theorem 3.3.4]{Ber-Jun-Mil}, every topological semigroup (group) $G$  has a universal totally disconnected semigroup compactification, which by Proposition \ref{nuB X Prop} (v) (Theorem \ref{First TD Semigroup Cpctn Thm}) is equivalent to  $(\nub G, \nub)$ for some Boolean subalgebra $\B$ of $\fB(G)$.  When $G$ is not a topological group, we have left the description of this Boolean subalgebra $\B$ to future work.
 \er

\subsection{The universal TD semitopological semigroup compactification}

Let $G$ be a topological group. For $B \in \fB(G)$, let 
$$\LO(B) = \LO_G(B) = \{ s \c B: s \in G\} \quad \text{ and } \quad \RO(B) = \RO_G(B) = \{  B \c s : s \in G\}$$ be the left and right $G$-orbits of $B$. Let 
$$\Bwap = \{ B \in \Bluc: \LO(B) \text{ is relatively }\sigfb-\text{compact}\}.$$

\br \rm  (a) Our goal is to identify $\B$ such that $(\nub G, \nub)$ is the universal semitopological semigroup compactification of $G$.  By \cite[Corollary 1.4.4]{Ber-Jun-Mil}, every semitopological semigroup compactification of $G$ has the joint continuity property, so $\B$ will necessarily be contained in $\Bluc$ by Theorem \ref{Universality of nuLUC Thm}. Equivalently, we could have required that sets in $\Bwap$ belong to $\Blmc$; indeed readers wishing to avoid using \cite[Corollary 1.4.4]{Ber-Jun-Mil} can, throughout this subsection,  replace all references to $\Bluc$, $\sigluc$, etc. with $\Blmc$, $\siglmc$, etc. 

\smallskip 

\noindent (b) For $B \in \Bluc$, $\LO(B) \subseteq \Bluc$, so in our definition of $\Bwap$ we could have required that $\LO(B)$ be relatively $\sigluc$-compact and, arguing as below, obtained the same compactification of $G$ and, therefore, the same collection of sets $\Bwap$. We chose the above definition because we do not know if $\Bluc$ is $\sigfb$-closed in $\fB(G)$, so  our definition is formally weaker than the alternative. \er 

\blem \label{Bwap G-invariant Lemma}  The collection of sets $\Bwap$ is a $G$-invariant Boolean subalgebra of $\fB(G)$. Hence, the totally disconnected compactification $(\nuwap G, \nuwap) \in \K_0(G) $ exists. 
\elem 

\begin{proof}  Since $\LO(\emptyset) = \{\emptyset\}$ and $\LO(G)=\{G\}$, $\emptyset, G \in \Bwap$. Let $A, B \in \Bwap$, $t \in G$. Then $G \bs A, A\cup B, A \cap B, t \c A, A \c t \in \Bluc$. Let $(s_i)$ to be a net in $G$. Then there is a subnet $(s_{i_j})$ of $(s_i)$ and $C, D \in \fB(G)$ such that $s_{i_j} \c A \ra C$ in $\sigfb$ and --- passing to another subnet if necessary --- $s_{i_j} \c B \ra D$ in  $\sigfb$. By Lemma \ref{Weak and point topologies lemma} (iv), $s_{i_j} \c (G \bs A) = G \bs (s_{i_j}  \c A) \ra G \bs C$ in $\sigfb$ and   
$s_{i_j} \c (A \cup B) =s_{i_j}  \c A \cup s_{i_j} \c B  \ra C \cup D$ in  $\sigfb$, so $G\bs A, A \cup B \in \Bwap$; %$$s_{i_j} \c (A \cap B) =s_{i_j}  \c A \cap s_{i_j} \c B  \ra C \cap D$ $\sigfb$; and  $s_{i_j} \c (A \c t) =(s_{i_j}  \c A) \c t  \ra C \c t$ $\sigfb$. Hence similarly, $G \bs A, A \cup B,
 similarly, $A\cap B \in \Bwap$ and, by Lemma \ref{Weak and point topologies lemma}(vi), $A \c t \in \Bwap$. Since $\LO(t \c A) = \LO(A)$, $t \c A \in \Bwap$ as well. 
\end{proof} 

\bt \label{nuwap G universal semitop semigp cpctn Thm} Let $\B$ be  a Boolean subalgebra of $\fB(G)$. The following statements are equivalent:
  \bi \item[(i)] $(\nub G, \nub)$ is a semitopological semigroup compactification of $G$;   \item[(ii)] $\B$ is left introverted and contained in $\Bwap$;
 \item[(iii)]$\B$ is introverted and Arens regular.  
\ei 
In particular, $\Bwap$ is the largest introverted Arens regular Boolean subalgebra of $\fB(G)$ and $(\nuwap G, \nuwap)$ is the universal totally disconnected semitopological semigroup compactification of $G$. 
\et

\begin{proof}   Let $B \in \Bwap$. If $\V \in \nuwap G$, then $\V = \U \cap \Bwap$ for some $\U \in \nuluc G$ by Proposition \ref{nuB X Prop} (iv), and $\V \c B = \U \c B \in \Bluc$ since $\Bwap$ is $G$-invariant and $\Bluc$ is left introverted. We \it claim \rm that $T_B : \nuwap G \ra (\Bluc, \sigluc): \V \mapsto \V \c B$ is continuous. 

We first observe that $T_B$ is continuous with respect to the point topology $\rholuc$ on $\Bluc$. Indeed, if $\V_i \ra \V$ in $\nuwap G$ and $s \in \V \c B$, then $\V \in \lwap(B \c s)$, so $\V_i \in \lwap (B \c s)$, and therefore $s \in \V_i \c B$, eventually; similarly if $s \notin \V \c B$, $s \notin \V_i \c B$ eventually. By Lemma \ref{Weak and point topologies lemma} (iii), $\V_i \c B \ra \V \c B$ in $\rholuc$. 
It follows that $\LO_{\nuwap G} (B) : = \{\V \c B: \V \in \nuwap G\}$ is a $\rholuc$-compact subset of $\Bluc$. Since $\rholuc$ is Hausdorff, $\LO_{\nuwap G} (B)$ is $\rholuc$-closed and therefore, since $s \c B = \nuwap(s) \c B$, $\ov{\LO_G(B)}^{\rholuc} \subseteq \LO_{\nuwap G} (B)$. From the density of $\nuwap(G)$ in $\nuwap G$ and $\rholuc$-continuity of $T_B$, we obtain $\LO_{\nuwap G} (B) = T_B(\nuwap G)\subseteq \ov{T_B(\nuwap(G))}^{\rholuc} =\ov{\LO_G(B)}^{\rholuc}$, and therefore  $\ov{\LO_G(B)}^{\rholuc} = \LO_{\nuwap G} (B)$ is a $\rholuc$-compact subset of $\Bluc$. By Lemma \ref{Weak and point topologies lemma} (v),  $\ov{\LO_G(B)}^{\rho_{\fB}} = \LO_{\nuwap G} (B)$ is a $\rho_\fB$-compact subset of $\fB(G)$. 
However, $\LO_G(B)$ is relatively $\sigfb$-compact by hypothesis, and the Hausdorff topology $\rho_{\fB}$ is coarser than the Hausdorff topology $\sigfb$, so $\ov{\LO_G(B)}^{\sigluc}= \ov{\LO_G(B)}^{\rholuc} = \LO_{\nuwap G} (B)$ and the $\sigfb$ and $\rho_\fB$ topologies agree on this set; as $\LO_{\nuwap G} (B)\subseteq \Bluc$, $\sigluc$ and $\rholuc$ agree on $\LO_{\nuwap G} (B)$. Hence, $T_B: \nuwap G \ra (\Bluc, \sigluc)$ is continuous, as claimed. 

Let $\V \in \nuwap G$. Then $\V \c B \in \Bluc$. For $s \in G$, $R_{s^{-1}}(\V) \in \nuwap G$ by Lemma \ref{Weak and point topologies lemma} (vi), and it is easy to check that $s \c (\V \c B) = R_{s^{-1}}(\V)\c B$; hence, $\LO_G(\V \c B)  \subseteq \LO_{\nuwap G}(B)$, which is $\sigluc$-compact. Thus, $\V \c B \in \nuwap G$. It follows that $\Bwap$ is left introverted and for each $B \in \Bwap$, $T_B: \nuwap G \ra (\Bwap, \sigwap): \V \mapsto \V \c B$ is continuous. By Proposition \ref{General Arens regular Prop}, $\Bwap$ is introverted and Arens regular. 

Let $\B$ be a Boolean subalgebra of $\fB(G)$. Assume that $\B$ is  introverted and Arens regular,   and let $B \in \B$. By Proposition \ref{General Arens regular Prop} (ii) and Lemma \ref{Properties of operations Lemma} (iii), $\{ \V \c B: \V \in \nub G\}$ is $\sigma_\B$-compact --- therefore $\sigma_\fB$-compact --- and contains $\LO_G(B)$; hence, $B \in \Bwap$.  By Corollary \ref{Semitop semigp cpctn and Arens regularity Corollary}, we  have established the last statement of the theorem and that (iii) implies (ii). The  implications (ii) implies (i) and (i) implies (iii)  are consequences of Proposition \ref{nuB X Prop} (iv), Remark \ref{Cpctn Ppts Prsvd By Qtnts Remark} and  Corollary \ref{Semitop semigp cpctn and Arens regularity Corollary}. 
\end{proof}

\bc \label{nuwapG Corollary 1} Let $B \in \fB(G)$. The following statements are equivalent: 
 \bi 
\item[(i)] $B \in \Bwap$; 
\item[(ii)] $B \in \Bluc $ and $\LO_{\nuwap G}(B): = \{ \V \c B: \V \in \nuwap G \}$ is a $\sigfb$-compact subset of $\fB(G)$; 
\item[(iii)]  $B \in \B^{RUC}$ and $\RO_G(B)$ is a relatively $\sigfb$-compact subset of $\fB(G)$; 
\item[(iv)] $B \in \B^{RUC}$ and  $\RO_{\nuwap G}(B): = \{  B \c \U : \U \in \nuwap G \}$ is a $\sigfb$-compact subset of $\fB(G)$. 
\ei 
\ec 

\begin{proof}  By Theorem \ref{nuwap G universal semitop semigp cpctn Thm} and Proposition \ref{General Arens regular Prop} (also from the proof of Theorem \ref{nuwap G universal semitop semigp cpctn Thm}), statement (ii) follows from statement (i) and, since $\LO_G(B)$ is contained in $\LO_{\nuwap G}(B)$, statement (ii) implies statement (i). Also by Theorem \ref{nuwap G universal semitop semigp cpctn Thm} and Proposition \ref{General Arens regular Prop}, statement (i) implies statement (iv), which in turn implies statement (iii) because $\RO_G(B)$ is contained in $\RO_{\nuwap G}(B)$. Symmetric arguments to those found in the proof of Lemma \ref{Bwap G-invariant Lemma} and Theorem \ref{nuwap G universal semitop semigp cpctn Thm} show that $\B^{RWAP}$, the collection of sets  $B \in \B^{RUC}$ for which $\RO_G(B)$ is $\sigfb$-compact, is an introverted Arens regular Boolean subalgebra of $\fB(G)$. Hence $\B^{RWAP} \subseteq \Bwap$ by Theorem \ref{nuwap G universal semitop semigp cpctn Thm}.
\end{proof} 

\br  \rm For $B \in \B^{LUC} \cap \B^{RUC}$,   we have established the equivalence of relative $\sigfb$-compactness of $\LO(B)$ and $\RO(B)$ without help from the Grothendieck double-limit criteria, the Eberlein-Smulyan theorem, or any pre-existing theory of $WAP(G)$ or $G^{WAP}$, (much of which relies on these and other deep theorems). %In the final section of the paper,  we will provide some additional information about sets in $\Bwap$ by employing some of this theory. 
\er 

\bc \label{nuwap G Corollary 2}  The Arens regular introverted Boolean subalgebra $\Bwap$ of $\fB(G)$ is inversion invariant and $(\nuwap G, \nuwap)$ is the universal totally disconnected involutive semigroup compactification of $G$. 
\ec 

\begin{proof}  Let $B \in \Bwap$. Then $B^{-1} \in \B^{RUC}$ because $B \in \Bluc$. Since $\RO(B^{-1}) = \{ C^{-1}: C \in \LO(B)\} = \LO(B)^{-1}$,  $\RO(B^{-1})$ is relatively $\sigfb$-compact by Lemma \ref{Weak and point topologies lemma} (vii); by  Corollary \ref{nuwapG Corollary 1}, $B^{-1} \in \Bwap$. By Proposition \ref{Involution Prop}, $(\nuwap G, \nuwap)$ is an involutive semigroup compactification of $G$. Universality follows from Proposition \ref{Involutive semigroup compactifications Prop} and Theorem \ref{nuwap G universal semitop semigp cpctn Thm}. 
\end{proof}

  \br  \rm     If $G$ is a discrete group,  $\Bluc = \B^{RUC} = {\cal P}(G) = \fB(G)$ and $\beta G \equiv_G \nu_{\fB(G)} G$ by \cite[Proposition 4.7(g)]{Por-Woo}. Hence, $G^{LUC} \equiv_G \beta G \equiv_G \nuluc G$  and therefore, by the main result in \cite{Lau-Pym}, $Z_t(\nuluc G) = \nuluc(G)$. By Proposition \ref{Top centre vs left and right Arens product Prop}, $\Bluc$ is an example of a non-Arens regular introverted Boolean subalgebra of $\fB(G)$. In particular, $(\nuwap G, \nuwap) \precneqq_G (\nuluc G, \nuluc)$ in this case. 
 
 % \smallskip 
  
%  \noindent (c) In the next section we will see that $(\nuwap G, \nuwap)$ is a topologist's compactification of $G$ when $G$ is locally compact and totally disconnected. As with the proof of Proposition \ref{nuluc G topologist cpctn Prop}, this only requires showing that $H \in \Bwap$ for any open subgroup $H$ of $G$. 
  
 \er
 
 Let $\Omega =\Omega(G)$ denote the ``open coset ring of $G$", i.e., the  Boolean subalgebra of $\fB(G)$ generated by the open cosets in $G$. 
 
 \bp \label{Omega(G) contained in BWAP Prop} The open coset ring $\Omega$ is an inversion-invariant, Arens regular, introverted Boolean subalgebra of $\Bwap$. Thus, $(\nu_\Omega G, \nu_\Omega)$ is a totally disconnected involutive semitopological semigroup compactification of $G$ and $(\nu_\Omega G, \nu_\Omega)\preceq_G (\nuwap G, \nuwap)$. 
 \ep 
 
 \begin{proof} We first show that $\Bwap$ contains $\Omega$. Taking $H$ to be an open subgroup of $G$, it suffices to show that $H \in \Bwap$;  to this end we will establish that ${\cal C} := \LO(H) \cup \{\emptyset\}$ is $\sigma_\fB$-compact:  Let $(B_i)_{i \in I}$ be a net in ${\cal C}$. If $(B_i)$ contains a constant subnet, it contains a convergent subnet. Supposing that this is  not the case, we \it claim \rm that $B_i \ra \emptyset$ in $\sigfb$. Let $\V \in \nufb G$. Then $\emptyset \notin \V$. If $I_0:= \{ i \in I: B_i \in \V\}$ is cofinal in $I$, then, since the sets in ${\cal C}$ are disjoint, $B_i = B$ for some $B \in \LO(H)$ and all $i \in I_0$, contradicting our assumption that $(B_i)_{i \in I}$ has no constant subnet. Hence, $B_i \notin \V$ eventually, and the claim follows  by   Lemma \ref{Weak and point topologies lemma}(ii). Hence, ${\cal C}$ is $\sigfb$-compact. 
 
 Obviously $\Omega$ is $G$-invariant and inversion-invariant. To see that $\Omega$ is left introverted, let $B \in \Omega$, $\V \in \nu_\Omega G$. Taking $\U \in \nuwap G$ such that $\U \cap \Omega = \V$ (via Proposition  \ref{nuB X Prop} (iv)), $\V \c B = \U \c B \in \Bwap$; in particular  $\V \c B \in \fB(G)$. Supposing now that $B$ is a left coset of the open subgroup $H$ of $G$, by Lemma \ref{Properties of operations Lemma} (iv) it suffices to show that $\V \c B$ is an open coset. Let $r,s,t \in \V \c B$. Then each of $B \c r$, $B \c s$ and $B \c t$ are left cosets of $H$ belonging to $\V$, so we must have $B \c r = B \c s = B \c t \in \V$. Therefore, $B \c (rs^{-1}t) = (B \c s)(s^{-1}t) = B \c t \in \V$, whence $r s^{-1}t \in \V \c B$. It follows that $\V \c B$ is an open coset in $\fB(G)$, as needed. The proposition now follows from  Theorem \ref{nuwap G universal semitop semigp cpctn Thm}.
 \end{proof} 
 
 A topological group $G$ is \it non-Archimedean \rm if its open cosets form an open base for $G$. For example, every totally disconnected locally compact group is non-Archimedean \cite[Theorem 7.7]{HR}. For more information about non-Archimedean topological groups, including many examples, the reader is referred to \cite{Meg-Shl}. 
 
 \bc \label{Top Cpctn Cor}  Let $G$ be a non-Archimedean topological group. Then the semigroup compactifications $(\nu_\Omega G, \nu_\Omega)$, $(\nuwap G, \nuwap)$, $(\nuluc G, \nuluc)$ and $(\nulmc G, \nulmc)$ are all topologist's compactifications of $G$. 
 \ec  
 
 \begin{proof} In this case, $\Omega$ is an open base for $G$ and $\Omega \subseteq \Bwap \subseteq \Bluc \subseteq \Blmc$, so this is a consequence of Proposition \ref{nuB X Prop} (ii).  
 \end{proof} 
 
 \br \rm \label{Thanks to referee Remark} Our original version of Corollary \ref{Top Cpctn Cor} was stated only for totally disconnected locally compact groups and was proved in Section 5 with the aid of function-space arguments. We are grateful to the referee for suggesting the argument found in the first paragraph of the proof of Proposition \ref{Omega(G) contained in BWAP Prop}, and for bringing the class of non-Archimedean groups to our attention. 
 \er

\subsection{The universal TD topological group compactification of $G$}

Let $G$ be a topological group. %For a subgroup $H$ of $G$, $|G:H|$ is the index of $H$ in $G$ and $G/H$ denotes the set of \it right \rm cosets of $H$ in $G$. ; when appropriate, we will use $(G/H)_r$ and $(G/H)_l$ to respectively denote the sets of right and left cosets of $G$. 
 %For a subset $B$ of $G$, we use $T_B$ to denote a complete set of right coset representatives of $\BG$ without redundancies. Thus,  $|G / \BG| = |G:\BG| = |T_B|$.  
 Let 
$$\Bap =\{B \in \fB(G): |G:\BG| < \infty\}$$ 
where $|G:H|$ is the index of a subgroup $H$ in $G$. Let $\Nf(G)$ be the set of all finite-index  closed  normal subgroups of $G$.  

\bp \label{Bap TFAE Prop}  Let $B$ be a subset of $G$. The following statements are equivalent: 
\bi \item[(i)] $B \in \Bap$;
\item[(ii)] $B \in \fB(G)$ and $\RO(B)$ is finite; 
\item[(iii)] $B \in \fB(G)$ and $|G:G_B|$ is finite; 
\item[(iv)] $B \in \fB(G)$ and $\LO(B)$ is finite; 
\item[(v)] $B= \emptyset$ or $B$ is a (finite) union of cosets of some $N \in \Nf(G)$. 
\ei 
\ep 

\begin{proof}  (i) $\Leftrightarrow$ (ii): Let $T$ be a complete set of right coset representatives of $\BG$ without redundancies.  Since $B \c s = B \c t$ exactly when $\BG s= \BG t$, the map $T  \ra \RO(B) : t \mapsto B \c t$ is a bijection. 

\smallskip 

\noindent (i) $\Rightarrow$ (v):  Suppose that $B \in \Bap$ is nonempty. Since $|G:\BG|$ is finite, $\BG$ has only finitely many conjugates and therefore its normal core in $G$, $N=  \bigcap_{g \in G} g(\BG) g^{-1}$, is a finite-index normal subgroup of $G$. One can readily check that $\BG$ is closed because $B$ is so; hence $N \in \Nf(G)$. As $N$ is contained in $\BG$, $B = NB$, and we have statement (v). 

\smallskip 

\noindent (v) $\Rightarrow$ (i): Assume (v) and that $B$ is nonempty. Then $N \in \fB(G)$ and is contained in $\BG$, so $B \in \Bap$.

\smallskip 

\noindent A symmetric argument establishes the equivalence of statements (iii), (iv) and (v). 
\end{proof}

%We are grateful to the referee for suggesting the proof of the following lemma, which is shorter and less complicated than our original argument.  \blem \label{_BG and G_B index Lemma} For  a subset $B$ of $G$,  $|G:\BG| < \infty$ if and only if $|G: G_B| < \infty$. Thus, $\B \in \fB(G)$ belongs to $\Bap$ if and only if $|G:G_B|<\infty$, which happens exactly when $\LO(B)$ is finite. \elem \begin{proof} Let $N$ be the intersection of all $G$-conjugates of $\BG$, i.e., let $N$ be the ``normal core" of $\BG$ in $G$. Since $N$ is contained in $\BG$, $B = NB$ and, since $N$ is a normal subgroup of $G$, we have $BN= NB = B$; hence $N$ is contained in $G_B$. If $\BG$ has finite index, so too does its normal core $N$, and therefore $G_B$, which contains $N$, also has finite index. The reverse implication is similar.  \end{proof} 

\blem  \label{Bap introverted Lemma1} The collection  $\Bap$  is an Arens regular introverted inversion-invariant  Boolean subalgebra of $\Omega$, the open coset ring of $G$, (and therefore of $\Bwap$).
\elem 

\begin{proof}  Since $\Nf(G)$ is contained in $\Omega$, so too is $\Bap$ by Proposition \ref{Bap TFAE Prop}.  Hence,  $\Bwap$ contains $\Bap$, so by Theorem \ref{nuwap G universal semitop semigp cpctn Thm}  it suffices to show that $\Bap$ is an inversion-invariant left introverted Boolean subalgebra of $\fB(G)$. Clearly, $\emptyset, G \in \Bap$. Let $B, C \in \Bap$. Then $\BG \cap {_CG}$ is a finite-index subgroup of $G$ that is contained in both $_{B \cap C}G$ and $ _{B \cup C}G$. It follows that $B \cap C, B \cup C \in \Bap$. Also, $_{G \bs B}G = \BG$, so $G\bs B \in \Bap$ and  $B^{-1} \in \Bap$ by Proposition \ref{Bap TFAE Prop} (v). Thus, $\Bap$ is an inversion-invariant Boolean subalgebra of $\fB(G)$ and, since $\LO(t\c B) = \LO(B)$ and $\RO( B \c t) = \RO(B)$, Proposition \ref{Bap TFAE Prop} also tells us that $\Bap$ is $G$-invariant.  Finally, let $\V\in \nuap G$, $B \in \Bap$. Then $B \in \Bluc$, so $\V \c B \in \fB(G)$ by Lemma \ref{BLUC Lemma} (iii).  By Lemma \ref{Properties of operations Lemma} (ii) $\BG \leq _{\V \c B} G$, so  $\V \c B \in \Bap$. 
\end{proof} 

Thus, $(\nuap G, \nuap)$  is a semitopological semigroup compactification of $G$. For the next lemma,  note that because $\Bap \subseteq \Bluc$, $\BG \in \Bap$ whenever $B \in \Bap$. 

\blem  \label{Bap introverted Lemma2}  Let $\U, \V \in \nuap G$, $B \in \Bap$. Then $B \in \U \sq \V$ if and only if $B \c a \in \V$ and $\BG a \in \U$ for some $a \in G$.      \elem

\begin{proof}  Let $T$ be a complete (finite) set of right coset representatives of $\BG$ without redundancies.  Letting $T_\V = \{ a \in T : B \c a \in \V\}$, $\V \c B = \bigcupdot_{a \in T_\V} \BG a$ because $\RO(B)= \{ B \c a : a \in T\}$ and $B \c s = B \c a$ exactly when  $s \in \BG a$. Therefore, if $B \in \U \sq \V$, $\V \c B = \bigcupdot_{a \in T_\V} \BG a$ belongs to the (prime) ultrafilter   $\U$, so  $\BG a \in \U$ for some $a \in T_\V$. Conversely, suppose that $B \c a \in \V$ and $\BG a \in \U$. Then for each $s \in \BG a$, $B \c s = B \c a \in \V$, so $\V \c B$ contains $ \BG a$, which belongs to $\U$; therefore $\V \c B \in \U$, meaning $B \in \U \sq \V$.  
\end{proof}

\br   \label{First Bap Remark}  \rm %(a)  Since $B \c s = B \c t$ exactly when $\BG s= \BG t$, the map $T_B  \ra \RO(B) : t \mapsto B \c t$ is a bijection. Hence, $|G:\BG| = |\RO(B)|$. Similarly, $|G:G_B| = |\LO(B)|$. Thus, $B \in \Bap$ if and only if $\RO(B)$ is finite. \smallskip  \noindent (b) Since finite-index closed subgroups are open and, as one can check, $\BG$ is closed whenever $B$ is closed, it follows from (a) and the definitions of $\Bwap$ and $\Bluc$  that $\Bap \subseteq \Bwap$.\smallskip  \noindent (c) 
If $G$ is compact and totally disconnected, then $\Bap = \fB(G)$. Indeed, in this case $(\nu_{\fB(G)} G, \nu_{\fB(G)}) $ and $(\nuluc G, \nuluc)$ are topologist's compactifications of $G$ by Proposition \ref{nuB X Prop} (ii) and Corollary \ref{Top Cpctn Cor}. Hence, $(\nu_{\fB(G)} G, \nu_{\fB(G)}) \equiv_G (G, {\rm id}_G) \equiv_G (\nuluc G, \nuluc)$, so $\fB(G) = \Bluc$ by Proposition \ref{nuB X Prop} (iv). Thus, if $B \in \fB(G)$, $\BG$ is an open, therefore finite-index, subgroup of the compact group $G$; so, $B \in \Bap$. 

% $f_B: s \ra B \c s: G \ra \fB(G)$ is continuous with respect to the discrete topology on $\fB(G)$ by Lemma \ref{BLUC Lemma} (a)(iii); therefore $\RO(B) = f_B(G)$ is finite.  

\er

\bt \label{nuap G universal TD gp cpctn Thm} Let $\B$ be a Boolean subalgebra of $\fB(G)$. The following statements are equivalent: 
\bi \item[(i)]  $(\nub G, \nub)$ is a topological group compactification of $G$; 
\item[(ii)] $\B$ is left introverted and $|G:\BG|< \infty$  for each $B \in \B$; 
\item[(iii)]  $\B$ is an inversion-invariant Arens regular introverted Boolean subalgebra of $\Bap$.  \ei 
In particular,   $(\nuap G, \nuap)$ is the universal totally disconnected topological group compactification of $G$. 

\et 

\begin{proof}   By Lemma \ref{Bap introverted Lemma1}, $(\nuap G, \nuap)$ is an involutive semitopological semigroup compactification of $G$.  Let $\V \in \nuap G$, and suppose that $B \in \V^* \sq \V$. By Lemma \ref{Bap introverted Lemma2}, $B \c a \in \V$ and $\BG a \in \V^*$ for some $a \in G$, meaning $a^{-1} B \in \V$ and $a^{-1}\BG = (\BG a)^{-1} \in \V$. Hence, $a^{-1} B \cap a^{-1} \BG \neq \emptyset$ and therefore $B \cap \BG \neq \emptyset$. Taking $b \in  B \cap \BG$, $B = B \c b = b^{-1} B$, so $bB = B$; hence, there is some $b_1 \in B$ such that $b b_1=b$ and therefore $e_G = b_1 \in B$. Hence, $\V^* \sq \V \subseteq \nuap(e_G)$. Since $\V^* \sq \V$ and $\nuap(e_G)$ are  $\Bap$-ultrafilters, we obtain $\V^* \sq \V = \nuap(e_G)$, the identity of $(\nuap G, \sq)$. As $\V^{**} = \V$, $\V \sq \V^* = \nuap(e_G)$ as well. It follows that $\nuap G$ is a group with $\V^{-1} = \V^*$ for $\V \in \nuap G$. By \cite[Corollary 1.4.5]{Ber-Jun-Mil}, the compact, Hausdorff, semitopological group $\nuap G$ is a topological group.  

Letting $(\alpha G, \alpha)$ be a totally disconnected topological group compactification of $G$, by Proposition \ref{Involutive semigroup compactifications Prop}  there is  an introverted, Arens regular, inversion-invariant Boolean subalgebra $\B$ of $\fB(G)$ such that $(\alpha G, \alpha) \equiv_G (\nub G, \nub)$. Let $B \in \B$.  By Remark \ref{First Bap Remark}, $\RO(\lb(B)) = \{ \lb(B) \c \U : \U \in \nub G\}$ is finite.  Observe that \beq \label{nuap universal eqn}   \lb(B) \c \nub(a) = \lb (B \c a)  \qquad \text{for } a \in G.\eeq Indeed, $\V \in \lb(B \c a)$ if and only if $B \c a \in \V$, meaning $a \in \V \c B$; equivalently, $\V \c B \in \nub(a)$, or $B \in \nub(a) \sq \V$, which is the same  as $\nub(a)\sq \V \in \lb(B)$. Since $\nub G$ is a group, this is equivalent to $\V \in \nub(a)^{-1} \sq \lb(B) = \lb(B) \c \nub(a)$, and we have established (\ref{nuap universal eqn}). Consequently, $\{ \lb(B \c a): a \in G\}$ is also finite.  By the Stone representation theorem, $\lb: \B \ra \fB(\nub G)$ is a bijection, so $\RO(B) = \{ B \c a: a \in G\}$ is finite.  Hence,  $B \in \Bap$ by Proposition \ref{Bap TFAE Prop}. Since $\B \subseteq \Bap$,  we conclude that   $(\nub G, \nub) \preceq_G (\nuap G, \nuap)$. We have established the last statement of the theorem and that (i) implies (iii). That (iii) implies (ii) is trivial. If statement (ii) holds,  then $\Bap$ contains $\B$,  so $(\nub G, \nub)$ is a factor in $\SK(G)$ of $(\nuap G, \nuap)$ and statment (i) follows from Remark \ref{Cpctn Ppts Prsvd By Qtnts Remark}. 
\end{proof}

We turn now to identifying $(\nuap G, \nuap)$ with the profinite completion of $G$, and noting some consequences.  By Theorem \ref{nuap G universal TD gp cpctn Thm} this may already be clear to some readers, however we do not have a reference. We first recall the relevant definitions. 

 With respect to the ordering $M \preceq N$ if and only if $N \subseteq M$, $\Nf = \Nf(G)$ is a directed set, and for $M \preceq N$ the maps $\vpnm: G/N \ra G / M: sN \ra sM$ are (continuous, open) epimorphisms between finite, discrete groups satisfying $\vp_{NK} = \vp_{MK} \circ \vpnm$ whenever $K \preceq M \preceq N$. That is, $(G/N; \vpnm; M \preceq N \in \Nf)$ is a projective, or inverse, system of (Hausdorff topological) groups, e.g., see \cite[12.1.26]{Pal} and \cite[6.13]{HR}. Hence, the projective limit
$$\rho_f G:= \varprojlim_{N \in \Nf} G/N = \{ u = (u_N)_N \in \Pi_{N \in \Nf} G/N: \vpnm(u_N ) = u_M \text{ for } M \preceq N \in \Nf\}$$
is a closed subgroup of the compact totally disconnected group $\Pi_{N \in \Nf} G/N$ \cite[Theorem 6.14]{HR}; therefore, $\rho_f G$ is a compact totally disconnected group. Consider the map $\rho_f : G \ra \rho_f G: s \mapsto (q_N(s))_N$, where $q_N: G \ra G/N$ is the  quotient homomorphism. Letting $\vp_N: \rho_f G \ra G/N$ denote the projection homomorphism, $\vp_N \circ \rho_f = q_N$ is a continuous homomorphism for each $N \in \Nf$, so $\rho_f$ is a continuous homomorphism. Let $u \in \rho_f G$. For any finite subset ${\cal F}$ of $\Nf$, $N_{\cal F} = \bigcap_{M \in {\cal F}} M \in \Nf$; letting $s_{\cal F} \in G$ be such that $s_{\cal F} N_{\cal F} = u_{N_{\cal F}}$, $\vp_M(\rho_f(s_{\cal F})) = q_M(s_{\cal F}) = u_M = \vp_M(u)$ for each $M \in {\cal F}$. Hence, $(s_{\cal F})_{\cal F}$ is a net in $G$ such that $\rho_f(s_{\cal F}) \ra u$.  Thus, $(\rho_f G, \rho_f)$ is a totally disconnected group compactification of $G$, called the \it profinite completion \rm of the topological group $G$.

\bp  \label{nuap G = profinite completion Thm}  Let $G$ be a topological group. Then $(\nuap G, \nuap) \equiv_G(\rho_f G, \rho_f)$. 
   \ep 
   
\begin{proof}  By Theorem \ref{First TD Semigroup Cpctn Thm} it suffices to show that $\Bap = \B^{\rho_f}$, where $\B^{\rho_f} = \{ \rho_f^{-1} (D): D \in \fB(\rho_f G)\}$. From the above definition, the subcollection ${\cal S} = \{ \vp_N^{-1}(\{aN\}): a \in G, \ N \in \Nf\}$ of $\fB(\rho_f  G)$ is a subbasis for $\rho_f G$. This means that $\B_{\cal S}$, the collection of finite intersections of sets from ${\cal S}$, is a basis for $\rho_f G$ and, since each set $D$ in $\fB(\rho_f G)$ is compact, it is a finite union of sets from $\B_{\cal S}$. Thus, $\fB(\rho_f G)$ is the Boolean algebra of sets generated by ${\cal S}$, and it follows that $\B^{\rho_f}$ is the Boolean subalgebra of  $\fB(G)$ generated by $\{ \rho_f^{-1}(S): S \in {\cal S}\}$. But for $a \in G$ and $N \in \Nf$, $\rho_f^{-1}(\vp_N^{-1}(\{aN\})) = q_N^{-1}(\{aN\}) =  aN$, so $\B^{\rho_f}$ is the Boolean subalgebra of $\fB(G)$ generated by the cosets of  $N$ in $\Nf$.  By condition (v) of Proposition \ref{Bap TFAE Prop}, we conclude that $\B^{\rho_f} = \Bap$.  
\end{proof}

\br \label{Thanks to referee Remark 2}  \rm We are grateful to the referee for suggesting the main ideas used in Propositions \ref{Bap TFAE Prop} and \ref{nuap G = profinite completion Thm}, which simplified our original approach. 
\er 

A topological group $G$ is \it residually finite \rm if for each $s \in G$ with $s \neq e_G$, there is a continuous homomorphism $\vp$ of $G$ into a finite group $F$ such that $\vp(s) \neq e_F$; equivalently,  $\bigcap \Nf =\{e_G\}$. Discrete examples include free groups, finitely generated nilpotent groups, and finitely generated linear groups. Clearly, $G$ is residually finite exactly when $\rho_f$ is a monomorphism.  Moreover, $G$ is totally disconnected when it is residually finite, because $G_e$, the connected component of $e_G$, is contained in  every open subgroup of $G$ \cite[Theorem 7.4]{HR}, and therefore $G_e$ is contained in  $\bigcap \Nf$. Thus, one might expect  $(\rho_f G, \rho_f)$ to be a topologist's compactification of $G$  whenever $G$ is residually finite. However,  this is not the case.

\bc \label{nuap not topologist's cpctn Cor}  Let $G$ be a  locally compact  group. Then $(\rho_f G, \rho_f)$, equivalently $(\nuap G, \nuap)$, is a topologist's compactification of $G$ if and only if $G$ is totally disconnected and compact. 
\ec    

\begin{proof} If $G$ is compact and totally disconnected, then $G$ is residually finite by \cite[Theorem 7.7]{HR}. In this case, $\rho_f$ is a homeomorphism of $G$ onto $\rho_f G$.  If $(\rho_f G, \rho_f)$ is a topologist's compactification of $G$, then $G$ is residually finite, therefore totally disconnected. Suppose towards a contradiction that $G$ is not compact, and take $H$ to be a compact open subgroup of $G$ \cite[Theorem 7.7]{HR}. Let  ${\cal F}$ be a finite subset of $\Nf$.  Then, since $G$ is non-compact,  the closed finite-index normal subgroup $N_{\cal F} = \bigcap_{N \in {\cal F}} N$ is not contained in $H$, and we can therefore choose some  $s_{\cal F}$ in $ N_{\cal F} \bs H$. Since $\rho_f$ is one-to-one, $(\rho_f(s_{\cal F}))_{\cal F}$ is a net in $\rho_f(G) \bs \rho_f(H)$ and for each $N \in \Nf$, $\vp_N(\rho_f(s_{\cal F})) = s_{\cal F} N = e_G N = \vp_N(\rho_f(e_G))$ whenever ${\cal F} \supseteq \{ N \}$; therefore $\rho_f(s_{\cal F}) \ra \rho_f(e_G) \in \rho_f(H)$. Hence, $\rho_f(H)$ is not  open in $\rho_f(G)$.  
\end{proof} 

The following statement is  an immediate consequence of Lemma \ref{Bap introverted Lemma1} and   Corollaries  \ref{Top Cpctn Cor} and \ref{nuap not topologist's cpctn Cor}. 

\bc Let $G$ be a totally disconnected locally compact, non-compact group. Then $(\nuap G, \nuap) \precneqq_G (\nu_\Omega G, \nu_\Omega)$. 
\ec

In a standard way, Theorem \ref{nuap G universal TD gp cpctn Thm} tells us that $(\nuap G, \nuap)$ is the unique (up to equivalence) totally disconnected topological group compactification $(\gamma G, \gamma)$ of $G$ with the property that for any compact totally disconnected group $K$ and any continuous homomorphism $\vp: G \ra K$ there is a continuous (automatically homomorphic) map $\vp^\gamma: \gamma G \ra K$ such that $\vp^\gamma \circ \gamma = \vp$; moreover, $(\nuap G, \nuap)$ is the minimum semigroup compactification in $\SK(G)$ with this property. (One can similarly use Theorems \ref{Universal semigp cpctn Thm}, \ref{Universality of nuLUC Thm} and \ref{nuwap G universal semitop semigp cpctn Thm} to describe $(\nulmc G, \nulmc)$/$(\nuluc G, \nuluc)$/ $(\nuwap G, \nuwap)$ as the unique semigroup compactification/semigroup compactification with the joint continuity property/semitopological semigroup compactification to which continuous homomorphisms from $G$ into compact right topological semigroups/ditto/compact semitopological semigroup compactifications extend.) We also have the following consequence of Theorem  \ref{nuap G universal TD gp cpctn Thm} and Proposition \ref{nuap G = profinite completion Thm}:

\bc Let $G$ be a topological group. Then $(\nuap G, \nuap)$ is the unique totally disconnected topological group compactification $(\gamma G, \gamma)$  of $G$ with the property that  for any finite   group $F$ and any continuous homomorphism $\vp: G \ra F$, there is a continuous homomorphism $\vp^\gamma: \gamma G \ra F$ such that $\vp^\gamma \circ \gamma = \vp$; moreover, $(\nuap G, \nuap)$ is the minimum semigroup compactification in $\SK(G)$ with this property.
\ec 

\begin{proof}  Since finite (Hausdorff) groups are compact and totally disconnected, $(\nuap G, \nuap)$ has this property.  Suppose that $(\gamma G, \gamma)$ is any semigroup compactification of $G$ with this property. For each $N \in \Nf$, let $q_N^\gamma: \gamma G \ra G/N$ be a continuous homomorphism such that $q_N^\gamma \circ \gamma = q_N$. Then $Q^\gamma(x):= (q_N^\gamma(x))_{N \in \Nf}$ is a continuous homomorphism of $\gamma G$ into $\Pi_{N \in \Nf} G/N$ satisfying $Q^\gamma \circ \gamma = \rho_f$. Since $\rho_f G$ is closed in  $\Pi_{N \in \Nf} G/N$ and $\gamma(G)$ is dense in $\gamma G$, the continuous homomorphism $Q^\gamma$ maps $\gamma G$ into $\rho_f G$. Hence, $(\rho_f G, \rho_f) \preceq_G (\gamma G, \gamma)$. 
\end{proof} 

% https://iopscience.iop.org/article/10.1070/IM1969v003n04ABEH000807

% https://en.wikipedia.org/wiki/Residually_finite_group

%https://math.stackexchange.com/questions/2611512/profinite-completions-and-inverse-limits

% https://arxiv.org/pdf/1811.04394.pdf

% https://core.ac.uk/download/pdf/82122306.pdf

\bex  \rm The profinite completions $(\rho_f G, \rho_f)$ of many discrete groups have been identified,  so $(\nuap G, \nuap)$ has been identified in many instances. Some specific examples are: 

\smallskip 

\noindent (a)  $\rho_f \mathbb{Z} = \varprojlim_{m\in {\mathbb N}} {\mathbb Z}/ m{\mathbb Z}$. Alternatively,  $\rho_f \mathbb{Z} =  \Pi_p {\mathbb Z}_p$, where $p$ ranges over all primes and ${\mathbb Z}_p$ denotes the $p$-adic integers; here $\rho_f (n) = (n)_p$ for $n \in {\mathbb Z}$, e.g., see \cite[Proposition 1.5.2]{Wil}.  This identification of $(\nuap {\mathbb Z}, \nuap)$ was previously derived in \cite[Example 4.3]{Ili-Spr} from a known description of ${\mathbb Z}^{AP}$.  

\smallskip 

\noindent (b)  For $n \geq 3$,  $\rho_f SL_n({\mathbb Z}) = \Pi_p SL_n({\mathbb Z}_p)$, where $p$ ranges over all primes. 

\smallskip 
\noindent (c) If $G$ is any infinite group for which $\Nf(G) = \{G\}$, e.g., $SL_n(\R)$ and    $PSL_n(\R)$ for $n\geq 2$, $\rho_f G$ is the trivial group. 
\eex

\section{TD semigroup compactifications as Gelfand compactifications} 

Let $G$ be a topological group. Until now, we have developed all properties of totally disconnected compactifications without employing function spaces. However, every semigroup compactification of $G$ is, up to equivalence, a Gelfand compactification, i.e., the Gelfand spectrum of a left $m$-introverted $C^*$-subalgebra of $CB(G)$. In this section we describe the relationship between the totally disconnected semigroup compactifications $(\nub G, \nub)$ and Gelfand compactifications of $G$, and we identify $(\nu_* G, \nu_*)$ with quotients of the familiar compactifications $(G^*, \delta_*)$, where $*$ is $LMC$, $LUC$, $WAP$ or $AP$.  

For the convenience of the reader, we recall some basic definitions and facts from the Gelfand theory of commutative $C^*$-algebras, e.g., see \cite{Kan}. Taking $A$ to be a commutative, unital Banach algebra, its Gelfand spectrum, $\Delta(A)$, is the set of all non-zero algebra homomorphisms of $A$ into $\C$ endowed with the topology of pointwise convergence on $A$, (equivalently the relative weak$^*$-topology inherited from the dual space $A^*$), with respect to which $\Delta(A)$ is a compact Hausdorff space. The Gelfand representation $\Gamma: A \ra CB(\Delta(A)): a \mapsto \hat{a}$,  defined by $\hat{a}(\varphi) = \varphi(a)$ for $a \in A$ and $\varphi \in \Delta(A)$, is a norm-decreasing algebra homomorphism. In the case that $A$ is a unital commutative  $C^*$-algebra, $\Gamma$ is an isometric $*$-isomorphism of $A$ onto $CB(\Delta(A))$ \cite[Theorem 2.4.5]{Kan}. 

Let $X$ be a topological space, $\B$ a Boolean subalgebra of $\fB(X)$. Let $\SB(X)$ be the $\|\cdot\|_\infty$-closed linear span of $1_\B:= \{1_B: B \in \B\}$ in $CB(X)$, a $C^*$-subalgebra of $CB(X)$ since  $1_\B$ is self-adjoint and closed under multiplication. Letting $\DB = \Delta(\SB(X))$  and $\delta_\B: X \ra \DB$ the point evaluation map, the Gelfand transform $\Gamma : \SB(X) \ra CB(\DB): f \mapsto \hat{f}$ is  an isometric $*$-isomorphism with $\Gamma^{-1} = j_{\delta_\B}: CB(\DB) \ra \SB(X): g \mapsto g \circ \delta_\B$. Thus, $CB(\DB)$ is generated by its idempotents $1_{\fB(\DB)}$, and it follows that $1_{\fB(\DB)}$ separates the points of $\DB$. Consequently, $(\DB, \delta_\B)$ is a totally disconnected compactification of $X$. 

For $m \in \DB$, let $\U_m = \{ B \in \B : m(1_B) = 1\}$. For $B \in \B$, $m(1_B) \in \{0,1\}$, which can readily be used to show  that $\U_m$ is an ultrafilter on $\B$, i.e., $\U_m \in \nub X$. (For example, for $B \in \B$, $1 = m(1_X) = m (1_B) + m(1_{X \bs B})$, so either $B \in \U_m$ or $X\bs B \in \U_m$.) 

\bp \label{Delta_B-nub cpctn isomorphism Prop} The map $h: (\DB, \delta_\B) \ra (\nub X, \nub): m \mapsto \U_m$ is a compactification isomorphism. For $\U \in \nub X$, $m_\U:= h^{-1}(\U) \in \DB$ satisfies $m_\U(1_B) = 1$ if $B \in \U$ and $m_\U(1_B) = 0$ if $B \in \B \bs\U$.
\ep 

\begin{proof}  For $x \in X$, $B \in h(\db(x)) = \U_{\db(x)}$ if and only if $1 = \db(x)(1_B) = 1_B(x)$, which holds if and only if $B \in \nub(x)$; hence $h(\db(x)) = \nub(x)$. If $\U_m = \U_n$ for $m,n\in \DB$, then $m$ and $n$ agree on $1_\B$ and therefore on $\SB(X)$; hence $h$ is one-to-one. Suppose that $m_i \ra m$ in $\DB$ and $\U_m \in \lb(B)$ for some $B \in \B$. Then $m(1_B) =1 $, so $m_i(1_B) = 1$, equivalently $\U_{m_i} \in \lb(B)$, eventually. Thus, $h$ is continuous. It is now automatic that $h$ is a compactification isomorphism. 
\end{proof} 

\bc \label{Idempotents in SB(X) Corollary}  If $f \in \SB(X)$ is an idempotent, then $f = 1_B$ for some $B \in \B$.  
\ec 

\begin{proof}  By Proposition \ref{nuB X Prop} (v), $(\DB, \db) \equiv_X (\nu_{\B^\Delta} X, \nu_{\B^\Delta})$, where $B^\Delta = \{ \db^{-1}(D): D \in \fB(\DB)\}$. By Propositions \ref{nuB X Prop} (iv) and \ref{Delta_B-nub cpctn isomorphism Prop},  $\B^\Delta = \B$. Since $1_{\fB(\DB)}$ is the set of idempotents in $CB(\DB)$ and  $j_\db$ is an isometric $*$-isomorphism of $CB(\DB)$ onto  $\SB(X)$, $j_\db(1_{\fB (\DB)}) = 1_{\B^\Delta} = 1_\B$ is the set of idempotents in $\SB(X)$. 
\end{proof} 

\bp \label{Introversion vs m-introversion Prop} A Boolean subalgebra  $\B$ of $\fB(G)$ is left introverted if and only if $\SB(G)$ is left $m$-introverted; in this  case, the compactification isomorphism $h: (\DB, \db) \ra (\nub G, \nub): m \mapsto \U_m$ satisfies $\U_{m \sq n} = \U_m \sq \,  \U_n$. 
\ep 

\begin{proof} Suppose that $\SB(G)$ is left $m$-introverted and let $\U \in \nub G$, $B \in \B$, $s \in G$. Then $s \c 1_B = 1_{s \c B}$ and $1_B \c s = 1_{B \c s}$ are idempotents in $\SB(G)$, so $s \c B, B \c s \in \B$ by Corollary \ref{Idempotents in SB(X) Corollary}. Letting $m_\U = h^{-1}(\U) \in \DB$, $m_\U \c 1_B(s) = m_\U(1_{B \c s}) = 1$ if and only if $B \c s \in \U$, i.e., exactly when $s \in \U \c B$; otherwise $m_\U \c 1_B(s) = m_\U(1_{B \c s}) = 0$, which happens exactly when $s \notin \U \c B$. Hence $m_\U \c 1_B = 1_{\U \c B}$, an idempotent in $\SB(G)$; that $\U \c B\in \B$ follows from Corollary \ref{Idempotents in SB(X) Corollary}. Conversely, suppose that $\B$ is left introverted and let $m \in \DB$, $B \in \B$, $s \in G$. Then $s \c 1_B = 1_{s \c B}, 1_B \c s = 1_{B \c s} \in \SB(G)$ and one can check that $m \c 1_B = 1_{\U_m \c B}$, so $m \c 1_B \in \SB(G)$. By linearity and continuity of the maps $\SB(G) \ra \ell^\infty(G): f \mapsto s \c f, f \c s, m \c f$,  $\SB(G)$ is left $m$-introverted. %$s \c f, f \c s, m \c f \in \SB(G)$ for any $f \in \SB(G)$.  

In either case,  $h$ is a compactification isomorphism of semigroup compactifications, so it is automatically a semigroup homomorphism. However, this is also easy to see directly: $B \in \U_m \sq \, \U_n$ if and only if $\U_n \c B \in \U_m$, which happens exactly when $1= m(1_{\U_n \c B}) = m (n \c 1_B) = m \sq n(1_B)$; this happens if and only if $B \in \U_{m \sq n}$.  
\end{proof}

\br \label{sigma topologies comparison Remark} \rm  If $\B$ is a Boolean subalgebra of $\fB(G)$, $B_i \ra B$ in   $\sigma_\B$  if and only if $1_{B_i} \ra 1_B$ in $\sigma(\SB(G), \DB)$ by Proposition \ref{Delta_B-nub cpctn isomorphism Prop}. (Indeed,  
suppose that $B_i \ra B$ in $\sigma_\B$ and let $m \in \Delta_\B$. If $m(1_B) =1$, then $B \in \U_m$ so $B_i \in \U_m$, equivalently, $m(1_{B_i}) = 1$, eventually; similarly if  $m(1_B) = 0$, then $m(1_{B_i}) = 0$, eventually.  The other direction is equally straightforward.) Since any $m \in \Delta_\B$ extends to $M \in \Delta(CB(G)) = \beta G$, $\sigma(\SB(G), \DB)$ is the subspace topology inherited from $\sigma(CB(G), \beta G)$. In particular, $B_i \ra B$  in $\sigma_\fB$  if and only if $1_{B_i} \ra 1_B$ in  $\sigma(CB(G), \beta G)$. 
\er 

Let $f \in CB(G)$. Then $f \in LMC(G)$ (respectively, $f \in LUC(G)$) if $s \mapsto f \c s: G \ra CB(G)$ is $\sigma(CB(G), \beta G)$- (respectively, $\|\c\|_\infty$-) continuous; $f \in WAP(G)$ (respectively, $f \in AP(G)$) if its orbit under left --- equivalently right --- translation, $\LO(f) := \{ f \c s : s \in G\}$, is relatively $\sigma(CB(G), CB(G)^*)$- (respectively, $\| \c \|_\infty$-) compact.  The semigroup compactifications $G^{LMC}$,  $G^{LUC}$, $G^{WAP}$ and $G^{AP}$ are the associated Gelfand spectra of $LMC(G)$, $LUC(G)$,  $WAP(G)$ and  $AP(G)$,  each of which is an $m$-introverted $C^*$-subalgebra of $CB(G)$ \cite{Ber-Jun-Mil}. 

\bp \label{Characteristic functions of BLUC etc Prop} Let $B \in \fB(G)$. Then \bi
\item[(i)] $B \in \Blmc$ if and only if $1_B \in LMC(G)$; 
\item[(ii)] $B \in \Bluc$ if and only if $1_B \in LUC(G)$; 
\item[(iii)] $B \in \Bwap$ if and only if $1_B \in WAP(G)$; 
\item[(iv)] $B \in \Bap$ if and only if $1_B \in AP(G)$.
\ei
\ep 

\begin{proof} Statements (i) and (iii) follow from Remark \ref{sigma topologies comparison Remark} and, noting that a $\sigma(CB(G), \beta G)$-limit of idempotents is again an idempotent,    the characterization \cite[Theorem 4.2.3]{Ber-Jun-Mil} of $WAP(G)$.  Observe that the set of idempotents, ${\mathscr I}$,  in $CB(G)$  is $\| \c \|_\infty$-closed and, since  $\| 1_E - 1_F\|_\infty $  is 1 whenever $E$ and $F$ are sets with $E \neq F$, the relative $\| \c \|_\infty$-topology on ${\mathscr I}$ is the discrete topology. Hence, $\LO(1_B)$ is closed in $(CB(G), \| \c \|_\infty)$ and, since $1_B \c s = 1_{B \c s}$,  $\LO(1_B)$ with the relative  $\|\c \|_\infty$-topology is equivalent as a $G$-space  to $\RO(B)$ with the discrete topology.  Statement (iv) is now obvious and statement (ii) 
 follows from Lemma \ref{BLUC Lemma} (a)(ii). \end{proof}

In the next theorem, for a semigroup compactification $(\gamma G, \gamma)$ of $G$, $\S^\gamma(G) := CB(\gamma G) \circ \gamma$, a left $m$-introverted $C^*$-subalgebra of $CB(G)$ \cite[Theorem 3.1.7]{Ber-Jun-Mil} such that $(\gamma G, \gamma) \equiv_G (\Delta(\S^\gamma(G)), \delta_\gamma)$. 

\bt \label{Universal TD P-cpctns Thm}  Let $(\gamma G, \gamma) \in \SK(G)$ and let ${\cal P}$ be a property of right topological semigroups.
\bi \item[(i)] The collection  $\B^\gamma:= \{ \gamma^{-1}(D): D \in \fB(\gamma G)\} = \{ E : 1_E \in \S^\gamma(G)\}$ is a left introverted Boolean subalgebra of $\fB(G)$ and $(\nu_{\B^\gamma} G, \nu_{\B^\gamma})$ is the projective maximum of the set $\{(\alpha G, \alpha) \in \SK_0(G): (\alpha G, \alpha) \preceq_G (\gamma G, \gamma)\}$.

\item[(ii)]  Suppose that $(\gamma G, \gamma)$  is universal with respect to ${\cal P}$. If $(\alpha G, \alpha) \in \SK_0(G)$ has ${\cal P}$, then $(\alpha G, \alpha) \preceq_G(\nu_{\B^\gamma} G, \nu_{\B^\gamma})$. If $(\nu_{\B^\gamma} G, \nu_{\B^\gamma})$ has ${\cal P}$ --- e.g., if ${\cal P}$ is preserved by quotients in $\SK(G)$ --- then $(\nu_{\B^\gamma} G, \nu_{\B^\gamma})$ is the universal totally disconnected semigroup compactification of $G$ with ${\cal P}$. 
\ei 
\et 

\begin{proof} Since $j_\gamma: CB(\gamma G) \ra \S^\gamma(G): f \mapsto f \circ \gamma$ is an isometric $*$-isomorphism, $1_E \in \S^\gamma(G)$ if and only if $1_E = j_\gamma(1_D) = 1_{\gamma^{-1}(D)}$ for some $D \in \fB(\gamma G)$, and it is clear that $\B^\gamma$ is a Boolean subalgebra of $\fB(G)$. For $D \in \fB(\gamma G)$ and $s \in G$, $s \c \gamma^{-1}(D) = \gamma^{-1}(D \gamma(s^{-1}))$,  $ \gamma^{-1}(D)  \c s = \gamma^{-1}( \gamma(s^{-1})D) \in \B^\gamma$ because  the maps $x \mapsto \gamma(s^{-1}) x, x \gamma(s^{-1})$ are homeomorphisms on $\gamma G$. Hence, $\S^{\B^\gamma}(G)$, the closed linear span in $CB(G)$ of the idempotents in $\S^\gamma(G)$, is translation invariant. To see that $\S^{\B^\gamma}(G)$ is left $m$-introverted --- and therefore $\B^\gamma$ is left introverted by Proposition \ref{Introversion vs m-introversion Prop} --- it suffices to show that $m \c 1_E \in \S^{\B^\gamma}(G)$ for $m \in \Delta_{\B^\gamma}$ and $E \in \B^\gamma$. Taking $M \in \Delta(\S^\gamma(G))$ to be an extension of $m$ to $\S^\gamma(G)$, $m \c 1_E = M\c 1_E \in \S^\gamma (G)$ and $(M \c 1_E)^2 = M \c 1_E$; hence $m \c 1_E \in \S^{\B^\gamma}(G)$ as needed. Thus, $(\nu_{\B^\gamma} G, \nu_{\B^\gamma}) \in \SK_0(G)$ and $(\nu_{\B^\gamma} G, \nu_{\B^\gamma}) \equiv_G(\Delta_{B^\gamma}, \delta_{B^\gamma})\preceq_G(\gamma G, \gamma)$, since $\S^{\B^\gamma}(G) \subseteq \S^\gamma(G)$ \cite[Theorem 3.1.9]{Ber-Jun-Mil}. Suppose that $(\alpha G, \alpha) \in \SK_0(G)$ with $(\alpha G, \alpha) \preceq_G (\gamma G, \gamma)$. By Propositions \ref{nuB X Prop} (v) and \ref{Delta_B-nub cpctn isomorphism Prop}, $(\alpha G, \alpha) \equiv_G (\nu_{\B^\alpha} G, \nu_{\B^\alpha}) \equiv_G (\Delta_{\B^\alpha}, \delta_{\B^\alpha})$ and $\S^{\B^\alpha}(G) \subseteq \S^\gamma(G)$, so $\S^{\B^\alpha}(G) \subseteq \S^{\B^\gamma}(G)$; hence, $(\alpha G, \alpha) \preceq_G (\nu_{\B^\gamma} G, \nu_{\B^\gamma})$. This proves statement (i), and statement (ii) follows.  
\end{proof} 

Using the universal properties of the $LMC$, $LUC$, $WAP$ and $AP$ compactifications of  $G$ \cite{Ber-Jun-Mil},  from Remark \ref{Cpctn Ppts Prsvd By Qtnts Remark}, Proposition \ref{Characteristic functions of BLUC etc Prop} and Theorem \ref{Universal TD P-cpctns Thm} we  obtain some of the results from Section 4 (where our proofs demonstrated Stone compactification methods and were independent of known universal properties of $G^{LMC}$, $G^{LUC}$, $G^{WAP}$, and $G^{AP}$): 

\bc \label{Universal properties of nulmc G, etc via Gelfand} Let $G$ be a topological group. Then:
\bi
\item[(i)] $(\nulmc G, \nulmc)$ is the universal totally disconnected semigroup compactification of $G$;
\item[(ii)] $(\nuluc G, \nuluc)$ is the universal totally disconnected  semigroup compactification of $G$ with the joint continuity property;
\item[(iii)] $(\nuwap G, \nuwap)$ is the universal totally disconnected  semitopological semigroup compactification of $G$; and
\item[(iv)] $(\nuap G, \nuap)$ is the universal totally disconnected topological group compactification of $G$. \ei
\ec

Let $G$ be a locally compact group. Using a universal property of $(G^{AP}, \delta_{AP})$, Ilie and Spronk identified the universal totally disconnected topological group compactification of  $G$ with the quotient group $G^{AP}/(G^{AP})_e$, where $(G^{AP})_e$ is the connected component of the identity in the compact group $G^{AP}$ \cite{Ili-Spr}. In the same paper, these authors introduced the ``idempotent compactification" $(G^I, \lambda_I)$, a semitopological semigroup compactification of $G$: letting $B_I(G)$  denote the closed linear span of the idempotents in $B(G)$, the Fourier--Stieltjes algebra of $G$ \cite{Kan-Lau}, $G^I = \Delta (B_I(G))$ and $\lambda_I$ is the point evaluation map. Recall from Section 4.3 that $\Omega$ denotes  the open coset ring of $G$.

\bp  The involutive semitopological semigroup compactifications of $G$, $(G^I, \lambda_I)$ and $(\nu_\Omega G, \nu_\Omega)$, are equivalent, i.e.,   $(G^I, \lambda_I)\equiv_G (\nu_\Omega G, \nu_\Omega)$. \ep

\begin{proof}  As a closed translation-invariant subspace of $B(G)$, $B_I(G) = A_\lambda$, the Arsac-Fourier space associated to a continuous unitary representation $\lambda = \lambda_I$ on a Hilbert space $\H_\lambda$ \cite[Lemma 2.8.3]{Kan-Lau}; $\lambda_I$ is fully described in \cite{Ili-Spr}, but we do not need this description. Note that $B_I(G) = A_\lambda$ is a left introverted homogeneous $*$-subalgebra of $CB(G)$ in the sense of \cite[Section 2.2]{Spr-Sto}. (This is a special case of \cite[Theorem 2.18]{Spr-Sto}, however this is easy to see directly: letting $\xi*_\lambda \eta(s)$ denote the matrix coefficient function $\l \lambda(s) \xi | \eta\r$ for $\xi, \eta \in \H_\lambda$ and $s \in G$, $(\xi*_\lambda \eta) \c s = \xi*_\lambda  \lambda(s)^*\eta \in A_\lambda$ and for $m \in VN_\lambda = A_\lambda^*$, 
$$m \c (\xi*_\lambda \eta)(s) = \l m, \xi*_\lambda \lambda(s)^*\eta \r = \l m \xi | \lambda(s)^* \eta\r = m \xi *_\lambda \eta(s),$$
so $m\c (\xi*_\lambda \eta) \in A_\lambda = B_I(G)$.) By Proposition 2.3 and Theorem 2.4 of \cite{Spr-Sto}, the closure ${\cal E}_I(G)$ of $B_I(G)$ in $CB(G)$ is a left $m$-introverted $C^*$-subalgebra of $CB(G)$ and $(G^I, \lambda_I) \equiv_G(\Delta({\cal E}_I(G)), \delta_{{\cal E}_I})$. However, $\{1_E: E \in \Omega \}$ is the set of idempotents in $B(G)$ \cite{Hos}, so ${\cal E}_I(G) = \S^\Omega(G)$. By Proposition  \ref{Delta_B-nub cpctn isomorphism Prop}, $(G^I, \lambda_I) \equiv_G (\Delta_\Omega, \delta_\Omega) \equiv_G (\nu_\Omega G, \nu_\Omega)$.  \end{proof} 

\br \label{Universal CH cpctn Remark}  \rm Letting $\E(G)$ denote the closure of $B(G)$ in $CB(G)$, $(G^\E, \delta_\E):= (\Delta(\E(G)), \delta_\E)$, the ``Eberlein compactification", is the universal (CH)-compactification of $G$ \cite[Section 2.4 and Theorem 3.14]{Spr-Sto}.  By Theorem \ref{Universal TD P-cpctns Thm}, if $(\nu_{\B^\E} G, \nu_{\B^\E})$ is a (CH)-compactification where $\B^\E:= \{ E : 1_E \in \E(G)\}$, then it is the universal totally disconnected (CH)-compactification. However, since it is not obvious that being a (CH)-compactification is preserved by  quotients in $\SK(G)$, we do not know if $(\nu_{\B^\E} G, \nu_{\B^\E})$ actually is a (CH)-compactification. As noted in the proof of \cite[Theorem 3.14]{Spr-Sto}, $(G^I, \lambda_I) = (\nu_\Omega G, \nu_\Omega)$ is a (CH)-compactification. Although $\Omega$ is contained in  $\B^\E$, we do not know if they are equal. This motivates the first of the questions stated below. 
\er 

\noindent {\bf Questions} 

\medskip 

\noindent We highlight some remaining questions: 
\medskip 

\noindent 1.  Does $\Omega = \B^\E$? Is $(\nu_{\B^\E} G, \nu_{\B^\E})$ the universal totally disconnected (CH)-compactification? Is  $(\nu_\Omega G, \nu_\Omega)$ the universal totally disconnected (CH)-compactification? (An affirmative answer to the first question gives an affirmative answer to the the others; see Remark \ref{Universal CH cpctn Remark}.) 

\smallskip 

\noindent 2. When is $(\nu_* G, \nu_*)$ equivalent to  $(G^*, \delta_*)$ when $*$ is $LMC$, $LUC$, $WAP$ or $AP$; i.e., when are the latter compactifications totally disconnected, cf. \cite[Theorem 4.7(j)]{Por-Woo}? 

\smallskip

\noindent 3. Diagram (d) on page 212 of \cite{Ber-Jun-Mil} shows the containment relations between a variety of compactifications of topological groups, with most arrows known to be irreversible. As Stone compactifications, can totally disconnected analogues of some of the compactifications shown in this diagram, but   not considered herein,  be identified?  What does the diagram look like in this scenario? We showed that $\nu_{\fB(G)} G \succeq_G \nulmc G \succeq_G \nuluc G\succeq_G \nuwap G \succeq_G \nu_\Omega G \succeq_G \nuap G$, where the third and last inequalities can be strict even when $G$ is locally compact, and $\nulmc G \equiv_G \nuluc G$ when $G$ is locally compact. Are there examples of topological groups $G$  for which  $\nu_{\fB(G)} G \succneqq_G \nulmc G$ --- equivalently, $\fB(G)$ is not left introverted --- or  $\nulmc G  \succneqq_G \nuluc G$, cf. \cite[Examples 2.2.6 and 4.5.8]{Ber-Jun-Mil}? % or $\nuwap G \succneqq_G \nu_\Omega G$? 

\bigskip 

  \noindent {\bf Acknowledgements:} It is a pleasure to  thank to the referee whose   valuable suggestions resulted in significant improvements to this paper, e.g., see Remarks \ref{Thanks to referee Remark} and \ref{Thanks to referee Remark 2}. We are also grateful to the referee for directing us to several of the  references cited herein.

\bigskip 

  \noindent {\bf Funding:}  This work was supported by  the Natural Sciences and Engineering Research Council of Canada; Alexander Stephens was supported by a 2019 NSERC USRA and Ross Stokke was supported by NSERC grant RGPIN-2015-05044.

\noindent {\sc Department of Mathematics, University of Manitoba, Winnipeg, MB, Canada, R3T 2N2}; email address: {\tt  stephe47@myumanitoba.ca}

\bigskip 

\noindent {\sc Department of Mathematics and Statistics, University
of Winnipeg, 515 Portage Avenue, Winnipeg, MB, Canada, R3B 2E9}; email address: {\tt r.stokke@uwinnipeg.ca}

\end{document}